\documentclass[11pt]{article}          
\usepackage{graphicx}
\usepackage{hyperref}
\usepackage{amsmath}
\usepackage{fullpage}
\usepackage{amssymb}
\usepackage{amsfonts}
\usepackage{latexsym}

\newcommand {\R}   {{\rm I\!R}}

\newcommand{\hf}{\frac12}

\newcommand{\bfv}{ {\bf{v}}}

\newcommand{\bbR}{\mathbb{R}}

\newcommand{\bfA}{{\bf A}}
\newcommand{\bfD}{{\bf D}}

\newcommand{\bfP}{{\bf P}}
\newcommand{\bfG}{{\bf G}}
\newcommand{\bfR}{{\bf R}}
\newcommand{\bfJ}{{\bf J}}
\newcommand{\bfL}{{\bf L}}
\newcommand{\bfS}{{\bf S}}

\newcommand{\bfI}{{\bf I}}
\newcommand{\bfU}{{\bf U}}
\newcommand{\bfQ}{{\bf Q}}

\newcommand{\bfF}{{\bf F}}
\newcommand{\bfY}{{\bf Y}}
\newcommand{\bfX}{{\bf X}}

\newcommand{\bfd}{{\bf d}}
\newcommand{\bfe}{{\bf e}}

\newcommand{\bfs}{{\bf s}}
\newcommand{\bfm}{{\bf m}}
\newcommand{\bfr}{{\bf r}}

\newcommand{\bfw}{{\bf w}}

\newcommand{\bfq}{{\bf q}}

\newcommand{\bfx}{{\bf  x}}

\newcommand{\bfu}{{\bf u}}

\newcommand{\bft}{{\bf t}}

\newcommand{\bfp}{{\bf p}}

\newcommand{\Fcurly}{\mathcal{F}}
\newcommand{\Acurly}{\mathcal{A}}

\newcommand{\bfphi}{{\boldsymbol \Phi}}

\newcommand{\bfepsilon}{\boldsymbol \epsilon}

\newcommand{\grad}{{\boldsymbol \nabla}}

\graphicspath{{figures/}}

\date{July 27, 2017}

\newcommand*\samethanks[1][\value{footnote}]{\footnotemark[#1]}
\begin{document}

\title{A Multiscale Method for Model Order Reduction in PDE Parameter Estimation} 



\author{%
Samy Wu Fung\thanks{Department of Mathematics and Computer Science, Emory University,  Atlanta, GA, USA. \texttt{\{samy.wu,lruthotto\}@emory.edu}} \and Lars Ruthotto\samethanks[1]}







\maketitle

\begin{abstract}
Estimating parameters of Partial Differential Equations (PDEs) is of interest in a number of applications such as geophysical and medical imaging. Parameter estimation is commonly phrased as a PDE-constrained optimization problem that can be solved iteratively using gradient-based optimization. 
A computational bottleneck in such approaches is that the underlying PDEs need to be solved numerous times before the model is reconstructed with sufficient accuracy. One way to reduce this computational burden is by using Model Order Reduction (MOR) techniques such as the Multiscale Finite Volume Method (MSFV).
    
In this paper, we apply MSFV for solving high-dimensional parameter estimation problems.
Given a finite volume discretization of the PDE on a  fine mesh, the MSFV method reduces the problem size by computing a parameter-dependent projection onto a nested coarse mesh. A novelty in our work is the integration of MSFV into a PDE-constrained optimization framework, which updates the reduced space in each iteration. We present a computationally tractable way of explicitly differentiating the MOR solution that acknowledges the change of basis.
As we demonstrate in our numerical experiments, our method leads to computational savings for large-scale parameter estimation problems where iterative PDE solvers are necessary and offers potential for additional speed-ups through parallel implementation.
\\

\textbf{Keywords:} model order reduction, parameter estimation, PDE constrained optimization, optimal control, geophysical imaging.
\end{abstract}

\section{Introduction}
\label{intro}
PDE parameter estimation problems frequently arise in practical applications, e.g., in geophysical imaging~\cite{Wardhow1988,Parker1994,pratt1999,deymor,EpanomeritakisAkcelikGhattasBielak2008,schulz2011computational,haber2014computational}, medical imaging~\cite{Arridge1999,CheneyEtAl1999, ArridgeSchotland2009,LipponenEtAl2013,deSturlerEtAl2013}, and hydrology~\cite{KaipioSomersalo2006,OliverBook2008}. Commonly, parameter estimation is formulated as a PDE-constrained optimization problem and solved iteratively using gradient-based optimization. One challenge in these approaches is the large computational cost required to numerically solve the discretized PDEs. The costs of PDE solves are compounded for mainly two reasons: First, the PDEs are solved repeatedly for different parameters until the parameters are estimated with reasonable accuracy. Second, the PDE needs to be solved with sufficient accuracy  to provide high-quality reconstructions and reliable gradient information for efficient optimization.

One popular approach to reduce the computational burden in PDE-constrained optimization is to use Model Order Reduction (MOR); see, e.g.,~\cite{KunischVolkwein2008,FahlEtAl2000,NegriEtAl2013,HimpeOhlberger2015,gallivan2004model,grimme1996model,calo2016randomized,ghasemi2015fast} and the recent survey~\cite{BennerEtAl2014}. A key idea is to reduce the complexity of the PDE solves by projecting the controls (i.e., the parameter) and/or the states (i.e., the PDE solutions) onto lower-dimensional subspaces. Several techniques have been used recently to compute subspaces with good approximation properties, e.g., reduced basis methods~\cite{NegriEtAl2013}, moment matching~\cite{FengBenner2007}, empirical interpolation~\cite{BarraultEtAl2004,deSturlerEtAl2013}, Krylov subspace methods~\cite{OConnellEtAl2017}, and perhaps most commonly Proper Orthogonal Decomposition~(POD)~\cite{BuiThanhEtAl2008,KunischVolkwein2008,FahlEtAl2000,KaletaEtAl2010,LipponenEtAl2013}.
MOR techniques have also led to remarkable performance gains in sampling methods used in statistical formulations of PDE parameter estimation problems; see, e.g.,~\cite{GalballyEtAl2009,LiebermanEtAl2010,MartinEtAl2012,SpantiniEtAl2015}. 

In order to be effective for PDE-constrained optimization, the solutions obtained using MOR techniques need to accurately approximate solutions of the original problem for a sufficiently large set of parameters~\cite{BuiThanhEtAl2008}. For example, the success of POD-based methods relies on the availability of sufficiently many and well-distributed \emph{snapshots} that are solutions of the original problem for given parameters which are used to build an MOR basis. While sampling the parameter space is feasible in low-dimensional settings, the problem becomes more difficult or intractable for high-dimensional parameter spaces. Similar restrictions apply to interpolatory MOR techniques.
 
In this work, we use multiscale finite volume methods (MSFV)~\cite{ElfendievHou2009,lee2009adaptive,hajibeygi2011adaptive,hajibeygi2014compositional,jenny2009modeling,parramore2016multiscale,jiang2007multiscale,chung2016generalized,chung2015mixed,chung2011energy} to numerically solve PDE-constrained optimization problems with high-dimensional parameter spaces. 
 In MSFV the reduced subspace is built by computing an operator-induced interpolation between grid functions on a nested coarse mesh and solutions on the original fine mesh. Being operator-induced means that the basis vectors used to represent the subspace are informed by the parameters.  
The interpolation is computed by solving  fine mesh versions of the original problem separately for each coarse mesh block, which can be parallelized, and, in contrast to, e.g., POD,  can avoid solving the fine mesh problem altogether. 
The MSFV method bears similarity with other adaptive meshing techniques, e.g.,~\cite{LipnikovEtAl2004,HoreshHaber2011}, however, the homogenization applied to obtain the reduced problem is operator-dependent and can thus capture larger variety in the parameters.
MSFV methods have been employed successfully to reduce the costs of PDE solves in  porous media flow problems~\cite{HouWu1997,MacLachlanMoulton2006,KalchevEtAl2016,moyner2016multiscale,de2017multiscale,durlofsky2007adaptive,lunati2008multiscale,jenny2005adaptive,jenny2003multi} and more recently in electromagnetics~\cite{haber2014multiscale}.

In this work, we integrate MSFV methods into the projected Gauss-Newton approach presented in~\cite{haber2014multiscale,jInvRuthottoHaberTreister} and exemplarily show their potential using a geophysical imaging example.  
Rather than using a fixed reduced space throughout the iterative scheme, we use a computationally tractable way to adaptively update the basis within the derivative-based optimization process \cite{de2017multiscale}. 
In our discretize-then-optimize framework this requires explicit differentiation of the MSFV method. 
Derivatives of MSFV methods have also been approximately computed using interpolation~\cite{krogstad2011adjoint}, adjoint methods~\cite{fu2009multiscale, fu2010multiscale,fu2011multiscale}, and more recently using a general algebraic framework~\cite{de2017multiscale}.
In contrast to these works, we explicitly differentiate the solution of the discretized problem obtained with the MSFV solver with respect to the parameter of interest. 
Our work is most closely related to~\cite{de2017multiscale}, however, we use local sensitivity equations alleviating the need for automatic differentiation. 
We also demonstrate that a parallel implementation of the derivative operators is tractable.


The paper is structured as follows. In Sec.~\ref{sec2}, we introduce the general mathematical framework for the parameter estimation problem and the formulation of the Direct Current Resistivity (DCR) survey. In Sec.~\ref{sec3}, we give a general overview of model order reduction (MOR) in the parameter estimation framework as well as a brief overview of MSFV methods. In Sec.~\ref{sec4}, we show how to compute the derivatives of the reduced forward problem for the optimization scheme. In Sec.~\ref{sec5}, we demonstrate that the reduced forward problem yields results that are comparable to methods using the full mesh in a reduced runtime. In Sec.~\ref{sec6}, we conclude with a discussion of possible extensions of our method.

\section{Parameter Estimation}
\label{sec2}
In this section, we provide a brief description of the general framework for PDE parameter estimation and highlight some challenges associated with solving their numerical solution. Our discussion follows the discussions in~\cite{haber2014computational,jInvRuthottoHaberTreister}.  We also introduce the Direct Current Resistivity (DCR) model problem, which we use for our numerical experiments. The DCR problem is also known as the inverse conductivity problem and is thus closely related, e.g., to Electrical Impedance Tomography~\cite{Arridge1999,CheneyEtAl1999, ArridgeSchotland2009,LipponenEtAl2013}. 

\subsection{Forward Problem}
Let $\Omega \subset \mathbb{R}^d$ be the computational domain with $d=2,3$ and let  $m\colon \Omega \to \mathbb{R}$  a model function describing the parameter of interest. We consider measurement data of the form 
\begin{equation}
\label{contData}
  \mathbf{D}_{ij} = \Fcurly_{ij}(m) + \epsilon_{ij},
\end{equation} 
where the forward operator is given by 
\begin{equation*}
\Fcurly_{ij}(m) = (p_i, u_j(m)),
\end{equation*}
for each $i=1,\ldots,N_r, j=1,\ldots,N_s.$ Here,  $p_i\colon \Omega \to \mathbb{R}$ is the $i$th receiver function, $u_j\colon \Omega \to \mathbb{R}$ is the field induced by the $j$th source $q_j\colon \Omega \to \mathbb{R}$,  $\mathbf{D} \in \R^{N_r \times N_s}$ are the discrete measurements, and $N_r$ and $N_s$ are the number of receivers and sources, respectively. 
We assume that the measurements are given by the $L_2$ inner product between $p_i$ and $u_j$ plus some additive noise $\epsilon_{ij}$ that for simplicity is assumed to be Gaussian white noise. The field $u_j$ (also called the state) satisfies
\begin{equation}
\label{contPDE}
  \Acurly(m)u_j  = q_j,
\end{equation}
where $\Acurly(m)$ is the underlying PDE operator including boundary conditions. While this framework is generic and can be used for many geophysical and medical applications, we focus in the following on DC Resistivity survey.

The DC Resistivity survey is an exploration technique used in many geophysical applications \cite{mcgillivray1992forward}. It uses artificial sources to introduce a direct current into the ground, thereby creating an electric potential field. Measurements of these potential fields are then collected on the surface and are used to estimate the conductivity in the subsurface. 

In this case, the operator $\Acurly$ in~\eqref{contPDE} corresponds to the steady-state heterogeneous diffusion operator, and the governing equations are given by
\begin{align}
\label{eq:DCR}
  \nabla \cdot (\sigma(m(x)) \nabla u_j(x)) &= q_j(x), &x \in \Omega, \nonumber \\
  \nabla u_j(x) \cdot \vec{n}(x) &= 0, &x \in \partial\Omega, 
  \\
  u_j(x) \to 0 \text{ as } &\vert x \vert \to \infty &j=1,\ldots,N_s. \nonumber
\end{align}
Here, the parameter of interest is  the conductivity function $\sigma \colon \Omega \to \mathbb{R}$, $m$ is the model function that parametrizes $\sigma$, for example, $\sigma(m) = \text{exp}(m)$, and $u_j$ is the potential field induced by the $j$th source $q_j$. The goal is then to estimate the model $m$ given surface measurements. While~\eqref{eq:DCR} can be extended to tensor-valued conductivities, we limit the following discussion to the scalar case for simplicity.

\subsection{Discrete Forward and Inverse Problem}
We follow the discretize-optimize approach in \cite{haber2014computational,jInvRuthottoHaberTreister}. To this end, we discretize the domain $\Omega$ on a uniform mesh $\mathcal{M}_{h}$ consisting of $N_m$ cells and $N_n$ nodes. Here, $h>0$ is a parameter that denotes the discretization size. Also, let $\bfm \in \mathbb{R}^{N_m}$, $\bfu_j $, $\bfq_j$, $\bfp_i \in \bbR^{N_n}$, and $\bfepsilon_j \in \bbR^{N_r}$ be the discretizations of the model $m$, the field $u_j$, the receiver $p_i$, the source $q_j$, and the measurement noise $\epsilon_{ij}$ given in~\eqref{contData} and \eqref{contPDE} on $\mathcal{M}_{h}$, respectively. The discrete data measurements are given by
\begin{equation*}
  \bfD = \bfF(\bfm) + \bfepsilon,
\end{equation*}
with the forward problem
\begin{equation}
\label{forwardProbFine}
\bfF(\bfm) = \bfP^\top \bfA(\bfm)^{-1} \bfQ = \bfP^\top \bfU,
\end{equation}
where $\bfA \in \bbR^{N_n} \times \bbR^{N_n}$ is the discretized PDE operator, $\bfF$ is the discrete forward operator, $\bfP = [\bfp_1 \; \bfp_2 \; ... \; \bfp_{N_r}]$ is the receiver matrix that maps from the fields to the data, $\bfQ = [\bfq_1 \; \bfq_2 \; ... \; \bfq_{N_s}]$ are the sources, $\bfU = [\bfu_1 \; \bfu_2 \; ... \; \bfu_{N_s}]$ are the induced fields, and $\bfepsilon = [\bfepsilon_1 \; \bfepsilon_2 \; ... \; \bfepsilon_{N_s}]$ is the measurement noise.

In the optimal control literature, $\bfm$ corresponds to the discrete "control" and the fields, $\bfu_j$, are the discrete "states" \cite{gamkrelidze2013principles,sargent2000optimal}. We assume that the operator $\bfA$ is nonsingular and differentiable with respect to the model parameters. To this end, we incorporate the third boundary condition in~\eqref{eq:DCR} by fixing the value of $\bfu_j$ at one grid point.  This renders the fields $\bfu_j$ uniquely defined given the model $\bfm$  and due to the construction of $\bfA$  also differentiable with respect to the model. 

The parameter estimation problem then aims at estimating the underlying model given the data, $\bfD$, sources $\bfQ$, and receivers $\bfP$. To this end we eliminate the PDE constraint and consider the so-called reduced optimization problem
\begin{equation}
\label{optfull}
\begin{split}
    &\min_{\bfm} \quad \bfphi \left(\bfP^\top \bfA(\bfm)^{-1}\bfQ, \bfD \right) + R(\bfm),  \\ & \:\text{ s.t. } \quad \bfm_L \leq \bfm \leq \bfm_H,
    \end{split}
\end{equation}
where $\bfphi: \bbR^{N_r \times N_s} \times \bbR^{N_r \times N_s} \to \bbR$ is a misfit function, $R: \bbR^{N_m} \to \bbR$ is a regularization term, and $\bfm_L$ and $\bfm_H$ can be used to enforce physical bounds for the model parameters. 
For some common choices of misfit and regularization functions, we refer to~\cite{haber2014computational}, and without loss of generality focus on the sum-of-squared-differences misfit 
\begin{equation}
  \bfphi(\bfD_{\rm pred}, \bfD_{\rm obs}) = \hf \| \bfD_{\rm pred}-  \bfD_{\rm obs} \| \Vert^2_F,
\end{equation}
where $\|\cdot\|_F$ denotes the Frobenius norm. In addition, we consider Tikhonov regularization
\begin{equation}
  R(\bfm) = \frac{\alpha}{2} \|\bfL (\bfm - \bfm_{\rm ref})\|_2^2,
\end{equation}
where $\bfm_{\rm ref}$ is the reference model, and the operator $\bfL$ can be, e.g., a discrete derivative operator and $\alpha>0$ is a regularization parameter. 
Since the dimensions of $\bfm$ are typically very large, we approximately solve~\eqref{optfull} using the projected Gauss-Newton-PCG method described in~\cite{haber2014computational}.
\subsection{Sensitivity Computation}
A key ingredient of derivative-based methods for solving~\eqref{optfull} is the sensitivity matrix $\bfJ(\bfm)$, which characterizes how small changes in the model affect the measurements. Since similar techniques will be used to differentiate the multiscale basis, we briefly review the concept of sensitivity computations given in~\cite{haber2014computational}.
For any small perturbation $\delta \bfm \in \bbR^{N_m}$, the sensitivity matrix $\bfJ(\bfm)$ satisfies
\begin{equation*}
  \bfP^\top\bfu(\bfm + \delta \bfm) \approx \bfP^\top\bfu(\bfm) + \bfJ(\bfm)\delta \bfm + {\cal O} (\delta \Vert \bfm \Vert^2).
\end{equation*}
As a consequence, we have that the sensitivity matrix is given by
\begin{align}\label{eq:Jtemp} 
  \bfJ(\bfm) = \bfP^{\top} {\frac{\partial \bfu}{\partial \bfm}}.
\end{align}
To obtain the derivatives of the fields with respect to the model parameter, we use implicit differentiation. To simplify notation, we consider the case of a single source $\bfq$ and write the discretized PDE-constraint as
\begin{align*}
  \bfA(\bfm)\bfu = \bfq.
\end{align*}
Applying the product rule to differentiate both sides with respect to $\bfm$, we obtain
\begin{align*}
  \grad_{\bfm}(\bfA(\bfm)\bfu) + \bfA(\bfm){\frac {\partial \bfu}{\partial \bfm}} = 0.
\end{align*}
For nonsingular $\bfA(\bfm)$ this is equivalent to
\begin{align*}
  {\frac {\partial \bfu}{\partial \bfm}} = -\bfA(\bfm)^{-1}\left( \grad_{\bfm}(\bfA(\bfm)\bfu) \right).
\end{align*}
Finally, inserting this into ~\eqref{eq:Jtemp} gives the sensitivity matrix
\begin{align}
\label{eq:sens}
  \bfJ(\bfm) = - \bfP^\top \bfA(\bfm)^{-1}\left( \grad_{\bfm}(\bfA(\bfm)\bfu) \right).
\end{align}

As can be seen in \eqref{optfull} and \eqref{eq:sens}, each evaluation of $\bfphi$ and product with $\bfJ$ or its transpose requires one PDE solve per source. This renders solving~\eqref{optfull} very expensive, particularly for parameter estimation problems involving hundreds of thousands or even millions of sources; see, e.g,~\cite{haber2016solving,van2015penalty}. This observation motivates lowering the cost of the PDE solves through model order reduction techniques.

\section{Model Order Reduction in Parameter Estimation}
\label{sec3}
In this section, we introduce reduced order modeling \cite{bai2005reduced,chung2016adaptive} using Multiscale Finite Volume (MSFV) techniques in the context of parameter estimation problems. MSFV is only one way to lower the computational costs associated with PDE-constraints and we refer to~\cite{BennerEtAl2014} for a general overview.

\subsection{Model Order Reduction}
MOR can be applied both to the model, reducing the dimensionality of the nonlinear optimization problem~\eqref{optfull}, and the fields, reducing the dimensionality of the PDE-constraint. In the following, we assume that the model is represented efficiently, e.g., using a tensor or OcTree mesh as in \cite{haber2012adaptive}, and focus on the latter part.  A common theme in MOR is to project the PDEs onto a small, $k$-dimensional, subspace (where $k \ll n$) that is spanned by the basis 
\begin{align*}
\bfS_k = [\bfs_{1} \; \bfs_{2} \; \ldots \; \bfs_{k}] \in \mathbb{R}^{N_n \times k}.
\end{align*}
We can then replace the forward problem in (\ref{forwardProbFine}) with the following reduced forward approximation
\begin{align}
 \label{eq:fwdProbRed}
 \bfF_{k}(\bfm) = \bfP^{\top}\bfU_k(\bfm), 
 \end{align}
where 
\begin{equation*}
\bfU_k(\bfm) = \bfS_k\bfA_k(\bfm)^{-1}\bfS_k^{\top}\bfQ
\end{equation*}
is the reduced approximation of the fields and 
\begin{equation}
\label{eq:reducedPDE}
\bfA_k(\bfm) = \bfS_k(\bfm)^{\top}\bfA(\bfm)\bfS_k(\bfm)
\end{equation} 
is the reduced PDE. Inserting the reduced forward problem into~\eqref{optfull} yields the surrogate problem
\begin{equation}
\label{eq:optredfixed}
\begin{split}
&\min_{\bfm} \quad \bfphi_{k}  \left( \bfP^{\top}\bfS_k\bfA_k(\bfm)^{-1}\bfS_k^{\top}\bfQ, \bfD \right) + R(\bfm), \\
&\text{ s.t. } \quad \bfm_L \leq \bfm \leq \bfm_H.
\end{split}
\end{equation}
 Since the approximate fields $\bfU_k(\bfm)$ change for every iteration in the optimization scheme (\ref{eq:optredfixed}), a key challenge in solving the reduced optimization problem is finding a good basis $\bfS_k$ so that $\bfU_k(\bfm)$ is a good approximation of $\bfU(\bfm)$ for all possible $\bfm$.

Prominent examples for finding these bases are moment matching~\cite{FengBenner2007}, superposition of locally reduced models~\cite{LohmannEid2007}, matrix interpolation~\cite{PanzerEtAl2010,AmsallemFarhat2011}, Interpolatory Model Order Reduction \cite{BaurEtAl2011,de2015nonlinear}, Proper Orthogonal Decomposition \cite{HaasdonkOhlberger2008,WillcoxPeraire2002} and a Greedy procedure~\cite{ElmanLiao2013,LiebermanEtAl2010, BuiThanhEtAl2008}, all of which have been explored in parameter estimation. Most of these techniques use an offline phase in which the PDEs are solved for a large number of right hand sides and a basis is constructed from these solutions. While these methods have been shown to be effective for many problems, they are more difficult to apply for the problem at hand in which the dimensionality of $\bfm$ is typically in the order of millions and sampling the parameter space becomes intractable. 

\begin{figure}[t]
  \centering
  \includegraphics[width=\textwidth]{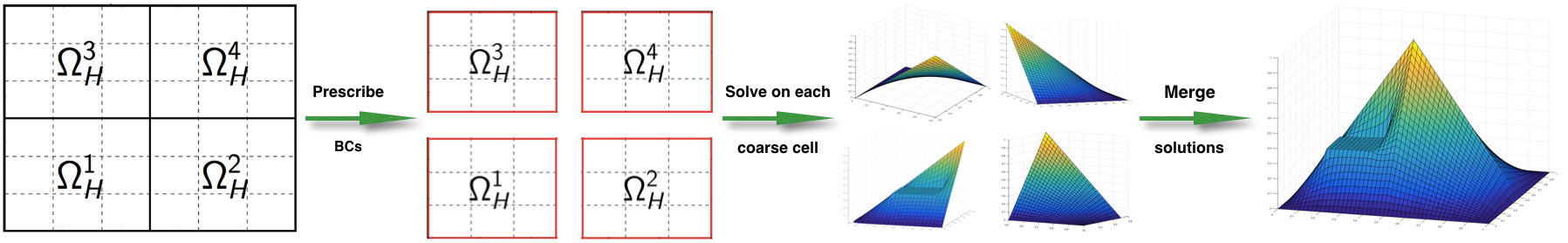}
  \caption{\textit{Outline of the construction of a multiscale basis function using a bilinear Piecewise Lagrangian Polynomial as boundary conditions. Here, $m(x)=10$ in $\Omega_H^1$, and $m(x)=1$ everywhere else. The adaptivity of the basis is evident by looking at its values on $\Omega_H^1.$}}
  \label{fig:fineToCoarse}
\end{figure}

To avoid the challenge of sampling the high-dimensional parameter space we propose using an adaptive basis 
\begin{equation*}
\bfS_k(\bfm) = [\bfs_{1}(\bfm) \; \bfs_{2}(\bfm) \; \ldots \; \bfs_{k}(\bfm)],
\end{equation*}
where now the projection basis depends on the model parameter $\bfm$. This leads to the optimization problem
\begin{align}
\label{eq:optredadapt}
&\min_{\bfm} \quad \bfphi_{k} \left( \bfP^{\top}\bfS_k(\bfm)\bfA_k(\bfm)^{-1}\bfS_k(\bfm)^{\top}\bfQ , \bfD \right)
+ R(\bfm), \nonumber \\ 
&\text{ s.t. } \quad \bfm_L \leq \bfm \leq \bfm_H.
\end{align} 
In the following subsection, we use the multiscale finite volume method \cite{ElfendievHou2009} to construct the basis $\bfS_k$. As we show in the next section, where we compute the sensitivities of  \eqref{eq:optredadapt} for MSFV, this leads to a computationally tractable smooth optimization scheme. While the adaptive basis approach is not limited to MSFV, computing derivatives might be more involved for other techniques.

\subsection{Multiscale FEM/FV Methods}
\label{sub:MultiscaleFEM}
We now introduce the basic steps of MSFV that are performed to construct the basis $\bfS_k$ to be used in~\eqref{eq:fwdProbRed}. Our discussion closely follows~\cite{ElfendievHou2009}. Each multiscale basis function is constructed in the following manner.
\begin{enumerate}
\item Partition the fine mesh $\mathcal{M}_{h}$ into a nested coarse mesh  $\mathcal{M}_{H} = \bigcup\limits_{j=1}^{N_c} \Omega_H^j$ where $\Omega_H^j$ is the $j$th coarse cell, and $N_c$ is the number of coarse cells. 
\item Choose a forcing term $q$, and prescribe values to a particular multiscale basis function $s: \bbR^{N_n} \to \bbR$ on the boundary of the coarse grid cells; see below for some common options. Denote these block-boundary values by $s_{\text{bc}}$.
\item Obtain the values of the multiscale basis function inside each coarse cell by solving the underlying PDEs on the local fine mesh using the prescribed boundary conditions for each coarse cell from step 2, i.e.,  
\begin{gather}
\label{eq:multiscaleBasis_step3}
  \Acurly(m) s = q, \:\: x \in \Omega_{H}^i, \:\: s = s_{\text{bc}}, \; x \in \partial\Omega_{H}^i,
\end{gather}
for each $i=1,...,N_c,$ where in the case of the DCR survey, $\Acurly$ corresponds to the diffusion operator, and $s$ is a particular multiscale basis function.
\end{enumerate}

The multiscale basis functions are thus obtained by solving the underlying PDEs \emph{locally} and \emph{independently} on each coarse mesh cell given specialized boundary conditions and/or forcing terms. We can therefore compute the bases in parallel. Each column in the projection basis $\bfS_k$ corresponds to one multiscale basis function; the number of columns in $\bfS_k$ therefore depends on the number of boundary conditions and/or forcing terms assigned. Furthermore,~\eqref{eq:multiscaleBasis_step3} shows that $\bfS_k$ is an operator-induced interpolation, thus allowing for the multiscale projection basis to adapt to the current optimal parameter in the optimization scheme. For a detailed analysis on the convergence of MSFEM/MSFV methods, we refer to \cite{hou1999convergence,efendiev2000convergence,chung2003convergence}.

The scheme is illustrated in Fig.~\ref{fig:fineToCoarse}. The specific choice of boundary conditions and/or forcing terms offers the flexibility to introduce application-specific prior knowledge into the construction.
We consider the following techniques:
\begin{itemize}
\renewcommand{\labelitemi}{$\bullet$}
  \item \textbf{Piecewise Lagrange Polynomials}: multiscale basis functions are constructed by setting $s_{\text{bc}}$ to be piecewise Lagrange polynomials on the coarse mesh, and solving \eqref{eq:multiscaleBasis_step3} with $q=0$. Due to the dependency of $\Acurly$ on $m$, the obtained basis function captures the local conductivity structure. For example, unless the PDE parameters are constant in the coarse mesh cell the obtained multiscale basis functions generally differ considerably from the generic FEM basis functions~\cite{ElfendievHou2009} .

  \item \textbf{Source Bases}: multiscale basis functions are constructed by solving \eqref{eq:multiscaleBasis_step3} with $q$ as the restricted sources from \eqref{contPDE} and setting $s_{\text{bc}}=0$ \cite{ElfendievHou2009,haber2014multiscale}.

  \item \textbf{Global Skeleton}: multiscale basis functions are constructed by first computing the fields $u_{\text{ref}}$ using a fixed reference parameter $m_{\text{ref}}$. The values of $u_{\text{ref}}$ on the coarse mesh (global skeleton) are then used as $s_{\text{bc}}$ in \eqref{eq:multiscaleBasis_step3} with $q=0$.

  \item \textbf{Local Bases}: multiscale basis functions are constructed by solving for the fields $u_{\text{ref}}$ using a fixed reference parameter $m_{\text{ref}}$. Subsequently  boundary conditions for each coarse mesh block are constructed. For example, on the $j$th cell, we apply a principal component analysis to the values of the fields $u_{\rm ref}$ on the boundary $\partial \Omega^j_H$ to identify the $r$ most important boundary conditions $u^j_{1}, u^j_{2}, \ldots,  u^j_{r}$. The associated multiscale basis functions are then obtained by solving~\eqref{eq:multiscaleBasis_step3} in $\Omega^j_H$ where $s^j_l = u_l^j$ on $\partial \Omega^j_H$ for each $l=1,2,\ldots,r$, and by keeping the values of the basis as zero in the rest of the domain. 
\end{itemize}
\section{Optimization with Multiscale FEM/FV Methods}
\label{sec4}
With the multiscale basis at hand, we now revisit the sensitivity computation of the reduced misfit in~\eqref{eq:optredadapt}. In this section, we derive the gradient, or sensitivities, of the reduced misfit and design an efficient mechanism for their implementation. A similar derivation has also been proposed recently in \cite{de2017multiscale}.

\subsection{Reduced Optimization}
We exemplarily compute the derivative of a data vector $\bfd$ obtained from a single source, $\bfq$, in~\eqref{eq:fwdProbRed}. The general case can then be obtained by decomposing the misfit function into a sum over the sources. 
Furthermore, we assume the sum-of-squared misfit
\begin{equation*}
  \bfphi(\bfm) = \hf\|\bfr_k(\bfm)\|^2, 
\end{equation*}
where the residual function is
\begin{equation*}
  \bfr_k(\bfm) = \bfP^{\top}\bfS_k(\bfm)\bfA_k(\bfm)^{-1}\bfS_k(\bfm)^{\top}\bfq - \bfd.
\end{equation*}
In this case, the
gradient is simply
\begin{equation*}
  \grad_{\bfm} \bfphi =  \bfJ_{k}(\bfm)^{\top} \ \bfr_k(\bfm),
\end{equation*}
where 
$\bfJ_{k}(\bfm)$ is the sensitivity (or Jacobian) of the \emph{reduced} forward model, that is,
\begin{align}
  \label{eq:Jk}
  \bfJ_{k}(\bfm)^\top = \grad_{\bfm} \left(
  \bfP^{\top}\bfS_k(\bfm)\bfA_k(\bfm)^{-1}\bfS_k(\bfm)^{\top}\bfq \right).
\end{align}
In the remainder of this section, we derive the sensitivity of the misfit for the case of a fixed basis in the reduced forward model.
We then give a detailed derivation of the (more complicated) gradient of the misfit when the basis is adapted to the model. Finally, we compare both results and provide an intuition about the difference of the gradients obtained for the fixed basis and adaptive basis.

\subsection{Optimization with Fixed Reduced Space} 
\label{sub:optimization_with_fixed_reduced_space}
For problems where $\bfS_k$ does not depend on $\bfm$, $\bfJ_{k}(\bfm)$ in \eqref{eq:Jk} can be computed by the
following sensitivity calculation. Writing 
\begin{equation*}
  \bfS_k^{\top}\bfA(\bfm)\bfS_k \bfu = \bfS_k^{\top} \bfq,  
\end{equation*}
and differentiating both sides with respect to $\bfm$, we obtain that
\begin{equation*}
  \bfS_k^{\top} \bfG + \bfS_k^{\top} \bfA(\bfm)\bfS_k\  {\frac {\partial \bfu}{\partial \bfm}} = 0, 
\end{equation*}
where the matrix $\bfG$ is obtained by differentiating $\bfA(\bfm)\bfS_k \bfu$ assuming $\bfu$ is constant, that is, 
\begin{equation}\label{eq:G}
  \bfG = \grad_{\bfm} (\bfA(\bfm)\bfS_k \bfu).
\end{equation}
This implies that
\begin{equation*}
  {\frac {\partial \bfu}{\partial \bfm}} = -(\bfS_k^{\top} \bfA(\bfm)\bfS_k)^{-1}\bfS_k^{\top} \bfG. 
\end{equation*}
Multiplying by the receiver matrix we see that the Jacobian of the reduced residual is
\begin{equation}
  \label{eq:senssimple}
  \bfJ_k(\bfm) = -\bfP^{\top}\bfS_k  (\bfS_k^{\top} \bfA(\bfm)\bfS_k)^{-1}\bfS_k^{\top} \bfG.
\end{equation}

\subsection{Optimization with Adaptive Reduced Space} 
\label{sub:optimization_with_adaptive_reduced_space}
Computing the sensitivity is more involved when using an adaptive reduced basis in which $\bfS_k$ depends on the model \cite{de2017multiscale}. In this case we need to differentiate the basis vectors,
$\bfs_{1}(\bfm),\ldots,\bfs_{k}(\bfm)$, with respect to the model $\bfm$.
This derivation is not standard and provided in detail below.  We also provide a description of our implementation that makes computing matrix-vector products with these derivatives tractable.

Similar to the previous section, the sensitivities for the general case are computed using implicit differentiation, i.e., differentiating both sides of
\begin{equation}\label{eq:redField}
  \bfS_k(\bfm)^{\top}\bfA(\bfm)\bfS_k(\bfm) \bfu = \bfS_k \left(\bfm \right)^{\top} \bfq.  
\end{equation}
For ease of presentation, we denote the operators that compute directional derivatives of $\bfS_k(\bfm)$ and its transpose by
\begin{equation}
  \label{eq:Y}
  \bfY_k(\bfv,\bfm) = \grad_{\bfm} \left(\bfS_k(\bfm)\bfv\right) 
\end{equation}
and
\begin{equation}
  \label{eq:X}
  \bfX_k(\bfw,\bfm) =  \grad_{\bfm} \left(\bfS_k(\bfm)^{\top}\bfw\right).  
\end{equation}
We also omit the dependency on $\bfm$ for brevity. Differentiating both sides of~\eqref{eq:redField} and using the notation in~\eqref{eq:Y} and \eqref{eq:X},  we obtain
\begin{align*}
  \bfX_k(\bfq) = &\bfX_k(\bfA \bfS_k \bfu) +  \bfS_k^{\top} \bfG \nonumber + \bfS_k^{\top}\bfA \bfY_k(\bfu) 
  \bfA_k {\frac {\partial \bfu}{\partial \bfm}},
\end{align*}
where $\bfG$ is as defined in~\eqref{eq:G}.
Assuming that the reduced discrete PDE is non-singular, this implies that
\begin{align*}
 {\frac {\partial \bfu}{\partial \bfm}} =  \bfA_k^{-1}\left(
\bfX_k(\bfq - \bfA \bfS_k \bfu) \nonumber - \bfS_k^{\top} (\bfG + \bfA \bfY_k(\bfu))\right).
\end{align*}
Multiplying by the receiver matrix, $\bfP^{\top}$, the current basis, $\bfS_k$, and applying the product rule yields the Jacobian of the reduced forward problem
\begin{equation}
  \label{eq:sensred}
  \bfJ_{k} = \bfP^{\top} \left( \bfY_k(\bft(\bfm)) + \bfS_k(\bfm) {\frac {\partial \bfu}{\partial \bfm}} \right ),
\end{equation}
where the coefficients of the fields with respect to the multiscale basis are denoted by
\begin{equation*}
  \bft(\bfm) = \left(\bfS_k(\bfm)^{\top}\bfA(\bfm)\bfS_k(\bfm)\right)^{-1}\bfS_k(\bfm)^{\top}\bfq.    
\end{equation*}
Clearly, there is a difference between the Jacobian obtained for a fixed basis in~\eqref{eq:senssimple} and the Jacobian obtained for the adaptive basis given in~\eqref{eq:sensred}. Note that even if $\bfS_k$ changes slowly with respect to $\bfm$, the term $\bfA \bfY_k(\bfu)$ may not be small as $\bfA$ is the discretization of a differential operator.
This might cause problems in the optimization when adapting the basis between iterates while using the Jacobian in \eqref{eq:senssimple}.
Indeed, if one ignores the dependence of $\bfS_k$ on $\bfm$ and uses the Jacobian in \eqref{eq:senssimple}, the error in the gradient can be rather large. 
\subsection{Illustrating the Error} 
\label{sub:erroranalysis}
\begin{figure}[t]
  \centering
  \includegraphics[height=2.7in, width=\textwidth]{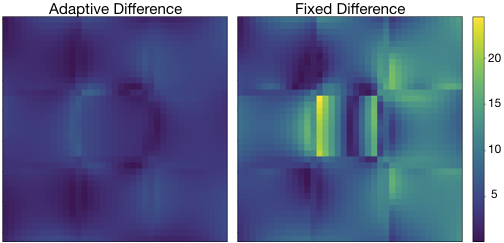}
  \caption{\textit{Difference plot of fine mesh and multiscale sensitivities. Here, we use one dipole source and the identity matrix as receivers. We obtain relative errors of 0.01 and 0.026 for the multiscale adaptive and fixed sensitivities, respectively. The color axis is chosen identically to allow for comparison.}}
  \label{fig:gradientError}
\end{figure}
Ignoring the dependency of $\bfS_k$ on the model $\bfm$ might not always lead to a large error in the gradient computation. 
Consider first the ideal case where the subspace $\bfS_k(\bfm)$ is chosen in such a way that
\begin{equation}\label{eq:MSexact}
\hspace{-0.8mm}  \bfP^{\top}\bfA(\bfm)^{-1}\bfq = 
  \bfP^{\top}\bfS_k(\bfm)\bfA_k(\bfm)^{-1}\bfS_k(\bfm)^{\top}\bfq 
\end{equation}
for every $\bfm$. In this case, it is clear that the dependence of $\bfS_k$ in $\bfm$ does not affect the quality of the multiscale solution.
We denote the error between the full and the reduced forward problem by
\begin{align*}
  \bfe_k(\bfm) = \bfP^{\top} \left( \bfA(\bfm)^{-1}\bfq -
  \bfS(\bfm)\bfA_k(\bfm)^{-1}\bfS(\bfm)^{\top}\bfq \right),
\end{align*}
and assume that for some $\epsilon>0$, the reduced forward problem satisfies the relaxed version of~\eqref{eq:MSexact}
\begin{equation*}
  \|\bfe_k(\bfm)\| \le \epsilon, \quad \forall \bfm.  
\end{equation*}
In this case, we have that for all $ \delta \bfm$,
\begin{align*}
  \|\bfe_k(\bfm + \delta \bfm) - \bfe_k(\bfm) \| &= \| (\bfJ_{k} - \bfJ) \delta \bfm \| + {\cal O} \|\delta \bfm\|^{2}  \\ &\le 2 \epsilon.
\end{align*}
Since the choice of $\delta \bfm$ is arbitrary, we obtain that the columns of $\bfJ$ and $\bfJ_k$ are at most
$2 \epsilon$ away.


Now, consider the solution of the optimization problem where the reduced model yields an accurate 
approximation to the desired function. In this case, if the change in $\bfm$ is not too large, ignoring the dependence of $\bfS_k$ on $\bfm$ may not lead to any serious complications. Nonetheless, if the subspace
approximation of the reduced model to the true model is not negligible, then ignoring it in the computation may lead to gross errors in the evaluation of the derivatives and to the lack of convergence of the optimization problem. 

As a small example, we compare the fine mesh sensitivities with the multiscale fixed and adaptive sensitivites for a test problem with mesh of size $36\times36\times12$, one dipole source, and the identity matrix as receivers. We choose the current model $\bfm$ so that $\sigma(\bfm)=0.1$ in the center to resemble a conductive block and $0.01$ in the rest of the domain. We compute the directional derivative for a constant perturbation $\delta \bfm=-0.001\:[1,\:...\:,1]^\top$. This results in relative errors of 0.01 for the adaptive sensitivities and 0.026 for the fixed sensitivities.  As to be expected, ignoring the dependence of $\bfS_k$ on $\bfm$ leads to big errors on the edge of the conductive block as can be seen in the difference plots in Fig.~\ref{fig:gradientError}.
\subsection{Computing $\bfY_k(\bfv, \bfm)$ and $\bfX_k(\bfw,\bfm)$}
The directional derivatives $\bfY_k(\bfv, \bfm)$ and $\bfX_k(\bfw,\bfm)$ are key components in the computation of the reduced adaptive sensitivities; see~\eqref{eq:sensred}. The computations of $\bfY_k(\bfv, \bfm)$ and $\bfX_k(\bfw,\bfm)$ require the multiscale basis $\bfS_k$ to be differentiable, which is ensured when using multiscale FEM/FV methods. 

Let the discretized version of~\eqref{eq:multiscaleBasis_step3} in a particular coarse cell be given by 
\begin{align*}
  \bfA (\bfm) \bfs_j(\bfm) = \bfq_j
\end{align*}
for $j=1,...,k,$ where $\bfs_j$ is the discretized version of the $j$th multiscale basis function. By construction of each multiscale basis function $\bfs_j, \: j=1,\ldots,k,$ this implies that 
\begin{align}
  \label{eq:DiscreteMSBasisPDE}
  \bfA (\bfm) \bfS_k(\bfm) \bfv =  \bfQ_k \bfv, 
\end{align}
where again, $\bfS_k(\bfm) = [\bfs_{1}(\bfm) \; \bfs_2(\bfm) \; \ldots \; \bfs_{k}(\bfm)]$, $\bfQ_k = [\bfq_{1} \; \bfq_2 \; \ldots \; \bfq_{k}]$, and  $\bfv \in \mathbb{R}^{k}$ is the vector corresponding to the directional derivative $\bfY_k(\bfv, \bfm)$. 

We exemplarily discuss the case of Dirichlet boundary conditions on the coarse block. Derivatives for local forcing terms can be computed along the same lines.
We reduce the linear system~\eqref{eq:DiscreteMSBasisPDE} in the following manner. For ease of presentation, denote \:$\bfS_k(\bfm) \bfv$ by \:$\bfx(\bfm) \in \mathbb{R}^{N_n}$. Let $N_I$ and $N_B$ be the number of interior and boundary nodes in the current cell, respectively, and define $\bfx_{I}(\bfm) \in \mathbb{R}^{N_I}$ and $\bfx_{B} \in \mathbb{R}^{N_B}$ as the corresponding values of $\bfx(\bfm)$ at the interior and boundary nodes of the chosen coarse cell (note that $\bfx_{B}$ is known and does not depend on the model). Similarly, let $\hat\bfq = \bfQ_k \bfv$, and define $\hat\bfq_{I} \in \mathbb{R}^{N_I}$ and $\hat\bfq_{B} \in \mathbb{R}^{N_B}$ as the entries of $\hat\bfq$ in the inner and boundary nodes, respectively.  We can then rewrite~\eqref{eq:DiscreteMSBasisPDE} in terms of the interior nodes as
\begin{align}
  \label{eq:reducedDiscreteMSBasisPDE}
  \bfA_{II}(\bfm) \bfx_I(\bfm) = \hat\bfq_I - \bfA_{IB}(\bfm)\bfx_{B},
\end{align}
where we partition the matrix $\bfA$ as follows
\begin{align*}
  \bfA = \begin{pmatrix} \bfA_{II} & \bfA_{IB} \\ \bfA_{BI} &  \bfA_{BB} \end{pmatrix}.
\end{align*}
Differentiating both sides of \eqref{eq:reducedDiscreteMSBasisPDE} with respect to $\bfm$, and using the product rule, we can isolate the term $\nabla_\bfm \bfx_{I}(\bfm)$ to obtain
\begin{align*}
  \nabla_{\bfm} \bfx_{I}(\bfm) = &- \bfA_{II}(\bfm)^{-1} (\nabla_{\bfm}(\bfA_{II}(\bfm)\bfx_{I}) + \nabla_{\bfm}(\bfA_{IB}(\bfm) \bfx_{B})),
\end{align*} 
which can be rewritten as 
\begin{align*}
  \nabla_{\bfm} \bfx(\bfm) = - \bfR_I \bfA_{II}(\bfm)^{-1} {\bfR_I}^\top (\nabla_{\bfm} (\bfA(\bfm)\bfx)),
\end{align*}
where $\bfR_I \in \mathbb{R}^{N_n \times N_I}$ is the basis for the inner nodes constructed by eliminating the columns of the identity matrix $\bfI \in \bbR^{N_n \times N_n}$ corresponding to the indices of the boundary nodes.
The above computation yields 
\begin{align*}\bfY_k(\bfv, \bfm) = \nabla_\bfm \bfS_k(\bfm)\bfv = \nabla_\bfm \bfx(\bfm)\end{align*} on a local coarse cell. This procedure is repeated for each coarse cell, and the solutions are merged together as it is done in the construction of $\bfS_k$. The construction of $\bfX_k(\bfw, \bfm)$ is along the same lines.

We like to emphasize that, similar to the construction of the multiscale basis $\bfS_k$, computing $\bfY_k(\bfv, \bfm)$ and $\bfX_k(\bfw, \bfm)$ can be done \emph{locally} and \emph{independently} on each coarse cell, leading to a parallelizable implementation. With $\bfS_k(\bfm)$, $\bfY_k(\bfv, \bfm)$, and $\bfX_k(\bfw, \bfm)$, we are thus able to solve the reduced adaptive optimization using gradient or Hessian based methods.

\section{Numerical Experiments}
\label{sec5}
In this section, we demonstrate the potential of the multiscale MOR parameter estimation method for a DC resistivity survey. We show that the multiscale inversion can reduce the time-to-solution compared to iterative PDE solvers, which are necessary for large-scale parameter estimation problems, with a moderate loss of reconstruction quality. We provide results for the fixed and adaptive multiscale basis; see Sec.~\ref{sub:optimization_with_fixed_reduced_space} and Sec.~\ref{sub:optimization_with_adaptive_reduced_space}, respectively.  We experiment on two test data: a block model that consists of two conductive blocks and that is homogeneous in the rest of the domain, and a 3D SEG/EAGE model of a salt reservoir described in \cite{aminzadeh1997}. We also perform a strong scaling test for the construction of $\bfS_k(\bfm)$, $\bfY_k(\bfv, \bfm)$, and $\bfX_k(\bfw, \bfm)$ to show the parallel efficiency of our current implementation.

Our multiscale framework is implemented as an extension to \texttt{jInv}~\cite{jInvRuthottoHaberTreister}, an open source package for PDE parameter estimation written in Julia \cite{bezanson2017julia}. For the discretization of the PDE operators, we use built-in methods in \texttt{jInv}, which are based on the mimetic finite volume method described in \cite{haber2014computational}. We use \texttt{jInv}'s methods for optimization, misfit functions, and regularizers. For brevity, we omit the term "multiscale" when referring to the multiscale adaptive and multiscale fixed inversions and refer to them as adaptive and fixed inversions instead. 
\subsection{Block Model Test Problem}\label{sub:ex2block}
We compare the reconstructions of the block model using the fine mesh inversion on a mesh of size $36\times36\times12$, and the multiscale inversion using fixed and adaptive bases for two different coarse meshes. The first one contains $9$ blocks of size $12^3$ each, and the second one contains $72$ blocks of size $6^3$ each. The test problem features $25$ sources and $1,369$ receivers located on the top surface. To construct the multiscale basis, we use boundary conditions obtained from Lagrange polynomials. We augment the basis by adding 25 basis functions of the global skeleton, and 93 and 400 local basis boundary conditions described in Sec.~\ref{sub:MultiscaleFEM} for the coarse meshes consisting of $12^3$ and $6^3$ blocks, respectively. Constructing these boundary conditions required $25$ fine mesh PDE solves in an offline phase which took about $0.24$ seconds. The overall number of basis functions is $k= 150$ and $k=572$ for the $12^3$ and $6^3$ respective coarsening strategies. 
To solve the fine mesh forward problem, we use MUMPS \cite{amestoy2000mumps} and a block CG method \cite{o1980block} with at most 100 iterations and stopping tolerance of $10^{-6}$. In order to make a fair comparison with the multiscale reconstructions, the stopping criteria for the block CG method are chosen such that the relative error of the block CG reconstruction has the same order of magnitude as the relative errors of the multiscale reconstructions; see Table \ref{tab:DC}. The reduced multiscale forward problems are all solved using MUMPS as we assume they are always small enough to be solved using a direct solver. For the inversions, we use 10 projected Gauss-Newton iterations with at most 15 CG iterations for each step. We add $1\%$ noise to the data, and enforce smoothness by using a diffusion regularizer with regularization parameter $\alpha=10^{-8}$.
\begin{figure}[!htb]
\newcommand{\rottext}[1]{\rotatebox{90}{\hbox to 45mm{\hss #1\hss}}}
    \begin{tabular}{ccc}
      &
      true model
      &
      full mesh reconstruction
      \\
      \\
      &
      \includegraphics[width=0.45\textwidth, height=1.6in]{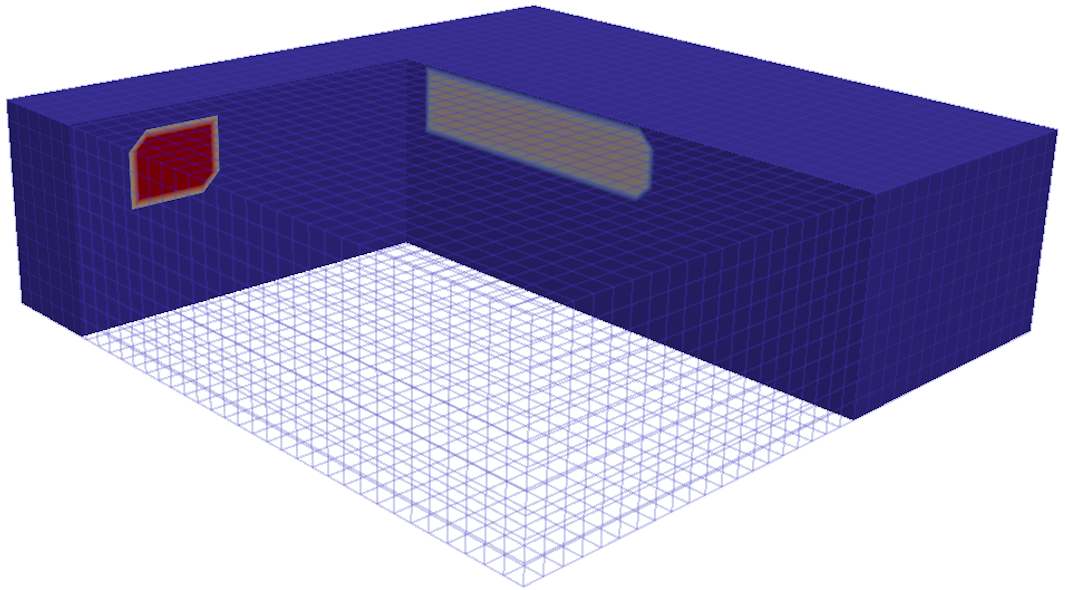}
       &
      \includegraphics[width=0.45\textwidth, height=1.6in]{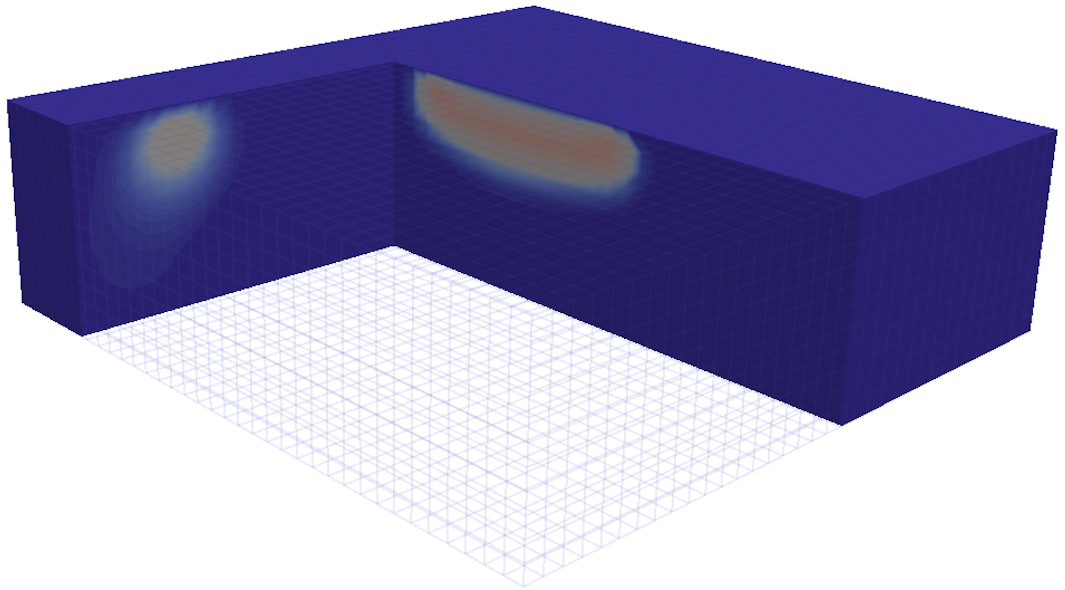}  
      \\ \hline
      \\
      &
      MS fixed reconstruction
      &
      MS adaptive reconstruction
      \\
      \rottext{$6\times6\times6$}
       &
      \includegraphics[width=0.45\textwidth, height=1.6in]{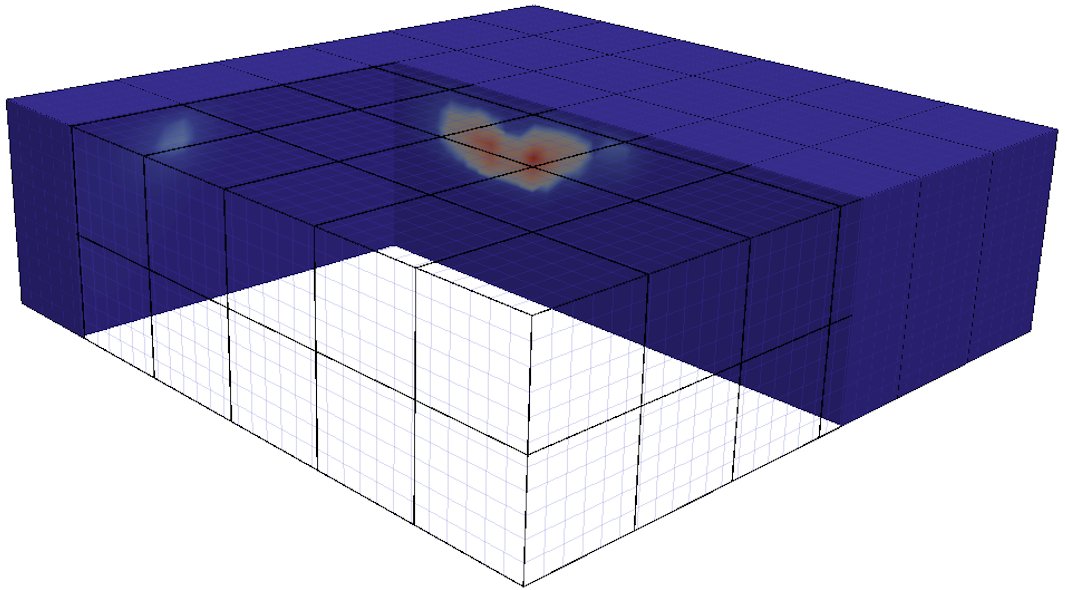}
       &
      \includegraphics[width=0.45\textwidth, height=1.6in]{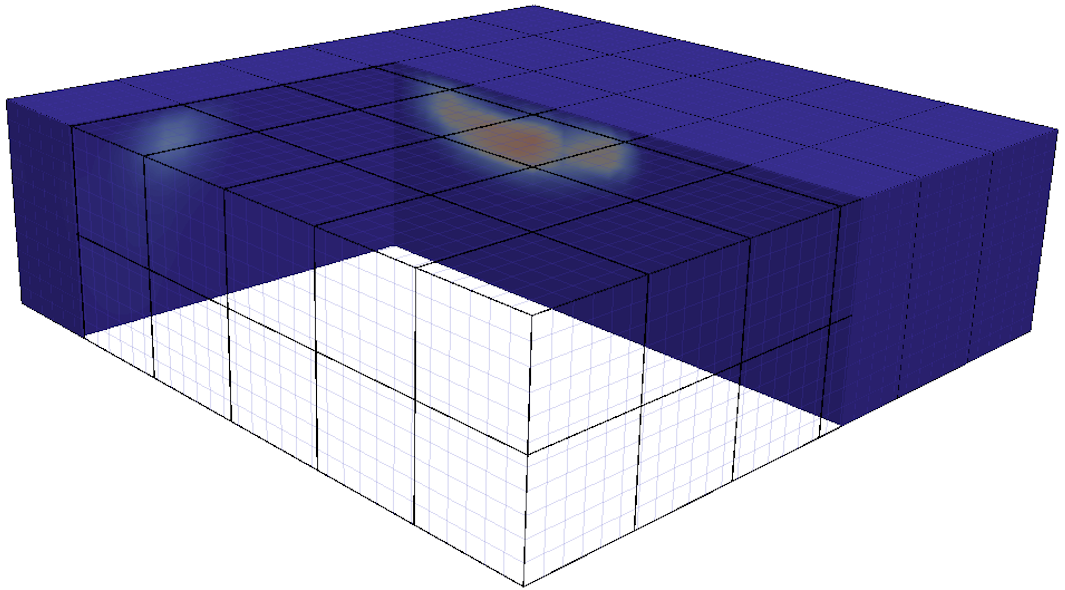} 
      \\ \hline
      \\
      &
      MS fixed reconstruction
      &
      MS adaptive reconstruction
      \\
      \rottext{$12\times12\times12$}
       &
      \includegraphics[width=0.45\textwidth, height=1.6in]{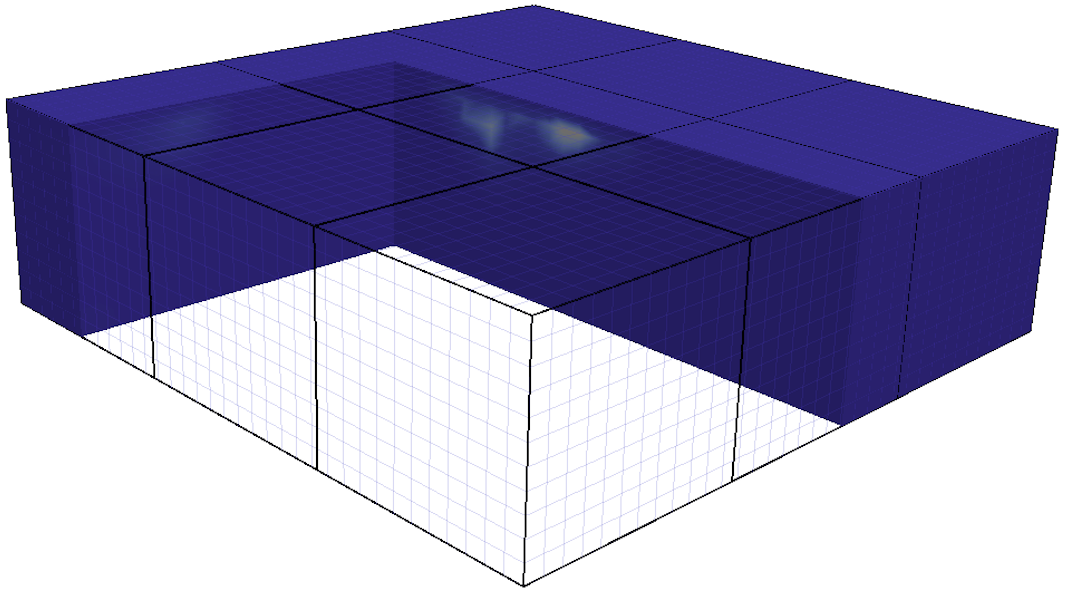}
       &
      \includegraphics[width=0.45\textwidth, height=1.6in]{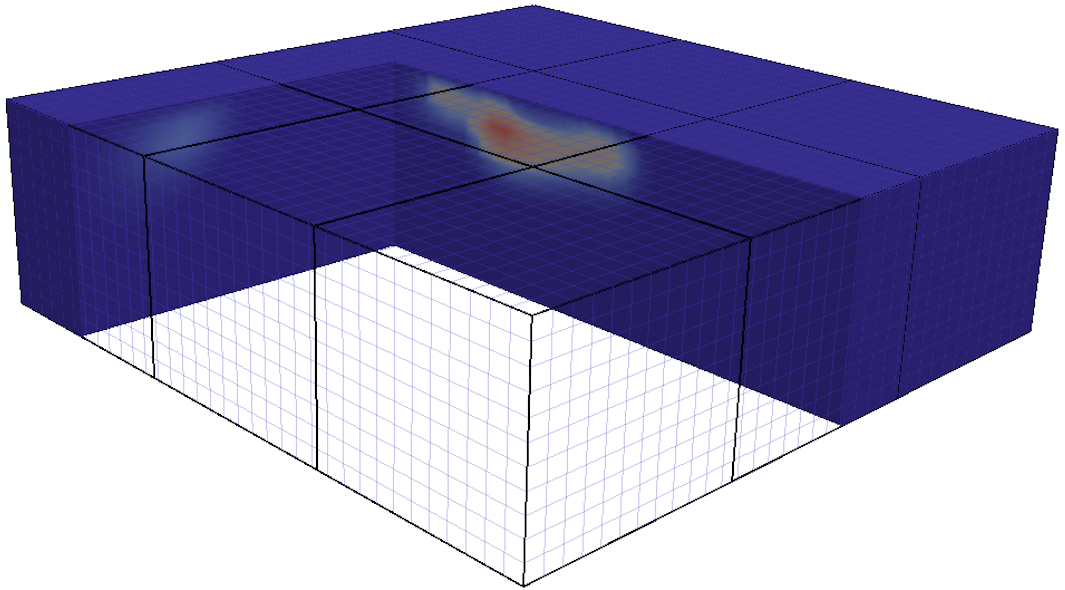} 
      \\
      \multicolumn{3}{c}{\includegraphics[width=0.7\textwidth, height=0.4in]{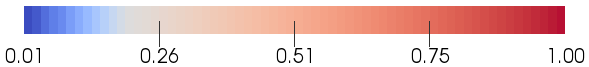}}
      \\ 
    \end{tabular}
    \caption{\textit{Model reconstructions for the block model test problem using the full mesh, and the adaptive and fixed multiscale inversions for two different coarsenings: $6^3$ and $12^3$ fine mesh cells per coarse cell. The figures were reproduced using the data visualization software Paraview \cite{ayachit2015paraview}.}}
  \label{fig:DCModels}
\end{figure}
\begin{table*}[t]
\centering
\caption{\textit{Relative errors and runtimes for the block model test problem in Sec.~\ref{sub:ex2block}. Here, $k$ corresponds to the total number of basis functions (and size of the PDE linear system), MS Fixed $6^3$ and MS Adaptive $6^3$ correspond to the multiscale inversions using coarse mesh blocks consisting of $6^3$ fine mesh cells, and MS Fixed $12^3$ and MS Adaptive $12^3$ correspond to the multiscale inversions using blocks of size $12^3$ fine mesh cells. The experiment is run on a standard Macbook air 2015 running macOS Sierra, with Intel core-i7 2.2 GHz CPU with 2 cores and 8 GB of RAM. The adaptive multiscale reconstruction is done in parallel using 2 processors. The fine mesh forward problem is solved using MUMPS and block CG. The reconstruction using MUMPS is considered as the baseline and is therefore assigned an error of $0$.}}
\label{tab:DC}
\begin{tabular}{|c|c|c|c|}
\hline
                     & $k$ & time (seconds) & relative error\\ \hline
fine mesh (MUMPS)    & 17,797               & 26.7     & 0\\ 
fine mesh (block CG) & 17,797               & 423.9   & 6.7e-2\\
MS fixed $6^3$       & 572              & 26.5     & 3.8e-2\\ 
MS adaptive $6^3$    & 572              & 306.6    & 2.1e-2\\ 
MS fixed $12^3$      & 150              & 25.4     & 4.9e-2\\ 
MS adaptive $12^3$   & 150              & 114.7    & 2.6e-2\\ \hline
\end{tabular}
\end{table*}
In Table~\ref{tab:DC}, we show results for the small block model test problem. The full mesh reconstruction requires solving a $17,797 \times 17,797$ linear system for the forward problem, whereas in the coarse meshes consisting of $6^3$ and $12^3$ cells, we project the PDEs down to $572 \times 572$ and $150 \times 150$ linear systems, respectively. We see that the adaptive inversions have a slower runtime than the fine mesh inversion when using MUMPS. Indeed, when the problem size is small enough to be solved using a direct solver such as MUMPS, there is no need to use MOR. However, in the typical setting where MOR is used, iterative solvers are more commonly used to solve large linear systems \cite{saad2003iterative}, and as we can see in Table~\ref{tab:DC}, both MS adaptive and MS fixed inversions are faster than the fine mesh inversion using the block CG algorithm as a linear solver. 

The adaptive inversion is in general slower than the fixed inversion since the projection bases must be rebuilt in every Gauss-Newton iteration and the sensitivity computations are more involved. However, as can be seen in the relative errors and in Fig.~\ref{fig:DCModels}, we obtain superior reconstructions using the adaptive inversions. The relative error of the adaptive reconstructions are about two times lower than the relative errors of the fixed reconstructions. This coincides with the images in Fig.~\ref{fig:DCModels}, where for the aggressive $12^3$ coarsening, the adaptive inversion manages to reconstruct the shape much better than the fixed reconstruction. For the more moderate $6^3$ coarsening, we project to a much larger subspace, leading to a fair reconstruction of the shape in both  fixed and adaptive inversions. However, the intensity values are more accurate in the adaptive reconstruction than in the fixed reconstruction.

\begin{table*}[h]
\centering
\caption{\textit{Relative errors and runtimes for the SEG test problem in Sec.~\ref{sub:SEG}. Here, $k$ corresponds to the total number of basis functions (and size of the reduced discretized PDE), MS Fixed $8^3$ and MS Adaptive $8^3$ correspond to the multiscale inversions using coarse mesh blocks consisting of $8^3$ fine mesh cells, and MS Fixed $16^3$ and MS Adaptive $16^3$ correspond to the multiscale inversions using blocks of size $16^3$ fine mesh cells. The experiment is run on a shared memory computer operating Ubuntu 14.04 with 2 × Intel
Xeon E5-2670 v3 2.3 GHz CPUs using 12 cores each, and a total of 128 GB
of RAM. Here,  Julia is installed and compiled using Intel Math Kernel Library
(MKL). The adaptive multiscale reconstruction is done in parallel using 16 processors. The fine mesh forward problem is solved using MUMPS and block CG. The reconstruction using MUMPS is considered as the baseline and is therefore assigned an error of $0$.}}
\label{tab:SEGmodels}
\begin{tabular}{|c|c|c|c|}
\hline
                     & $k$          & time (hours) & relative error\\ \hline

fine mesh (MUMPS)    & 139,425      & 0.46          & 0\\ 
fine mesh (block CG) & 139,425      & 5.98          & 3.2e-2\\
MS Fixed $8^3$       & 4415         & 0.25          & 7.0e-2\\ 
MS Adaptive $8^3$    & 4415         & 3.08          & 4.0e-2\\ 
MS Fixed $16^3$      & 1009         & 0.22          & 1.6e-1\\ 
MS Adaptive $16^3$   & 1009         & 1.97          & 1.1e-1\\ \hline
\end{tabular}
\end{table*}
\subsection{SEG/EAGE Test Problem}\label{sub:SEG}
\begin{figure}[!htb]
\newcommand{\rottext}[1]{\rotatebox{90}{\hbox to 45mm{\hss #1\hss}}}
    \begin{tabular}{ccc}
      &
      true model
      &
      full mesh reconstruction
      \\
      \\
      &
      \includegraphics[width=0.45\textwidth, height=1.6in]{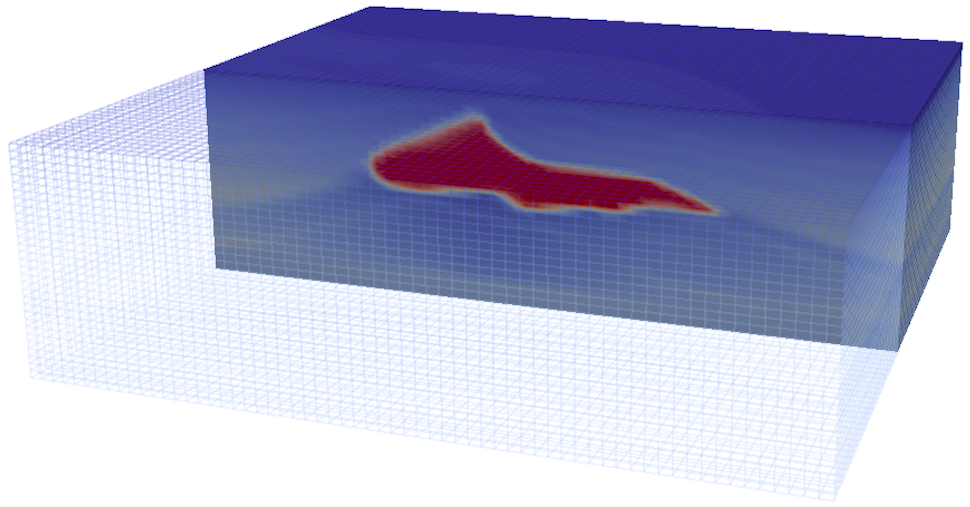}
       &
      \includegraphics[width=0.45\textwidth, height=1.6in]{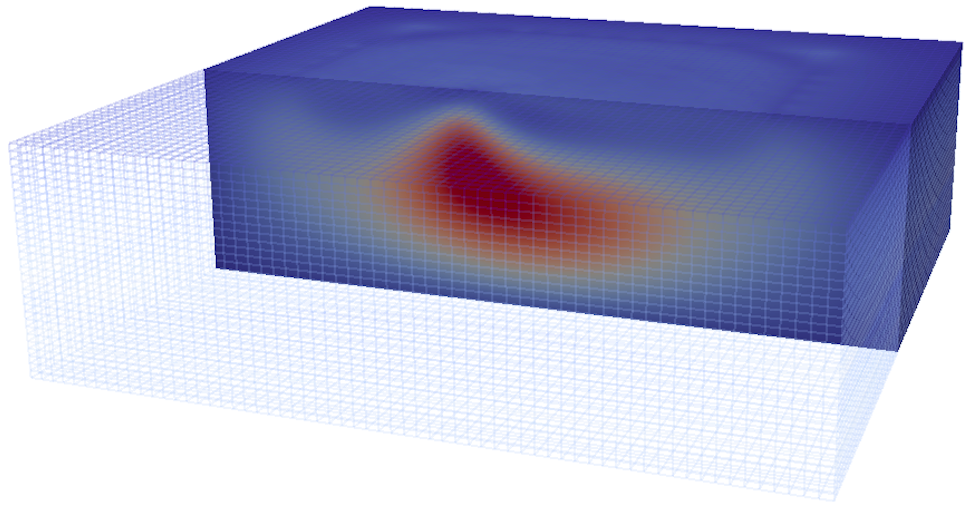}
      \\ \hline
      \\
      &
      MS fixed reconstruction
      &
      MS adaptive reconstruction
      \\
      \rottext{$8\times8\times8$}
       &
      \includegraphics[width=0.45\textwidth, height=1.6in]{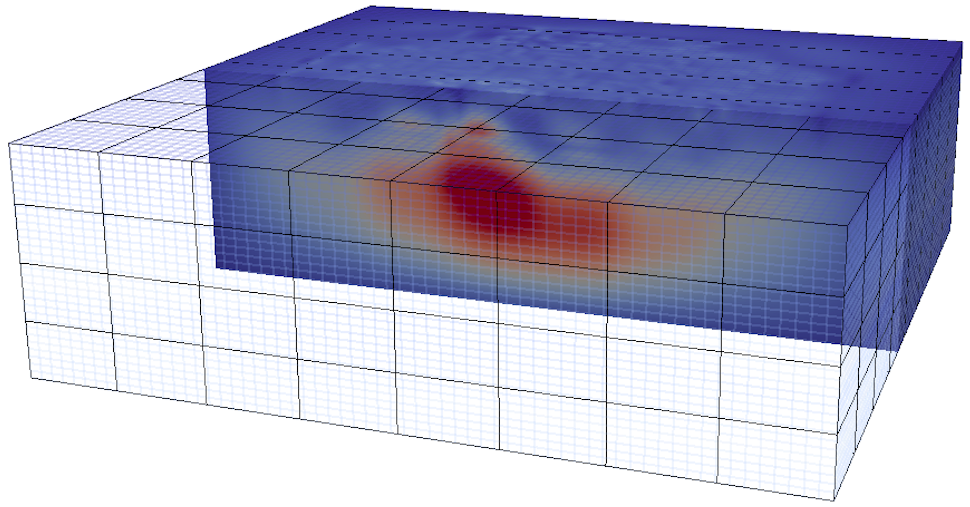}
       &
      \includegraphics[width=0.45\textwidth, height=1.6in]{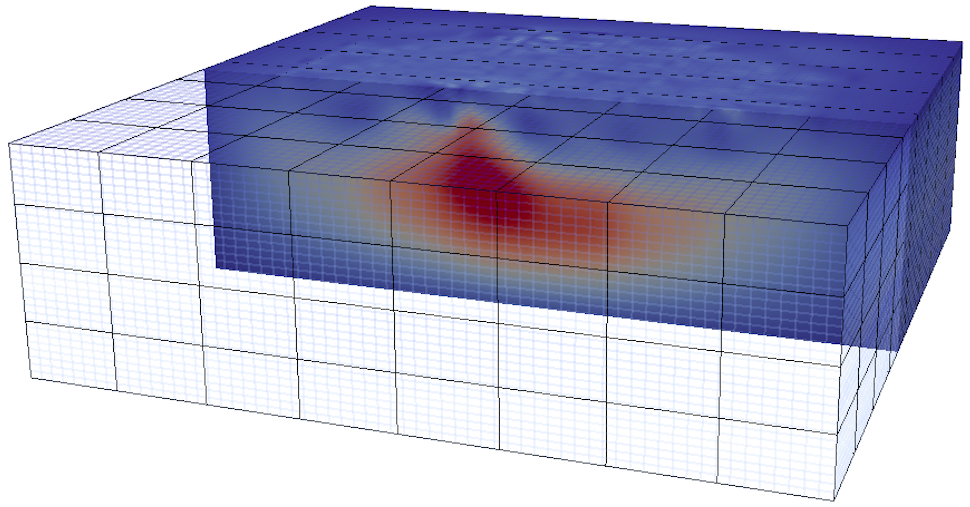}
      \\ \hline
      \\
      &
      MS fixed reconstruction
      &
      MS adaptive reconstruction
      \\
      \rottext{$16\times16\times16$}
       &
      \includegraphics[width=0.45\textwidth, height=1.6in]{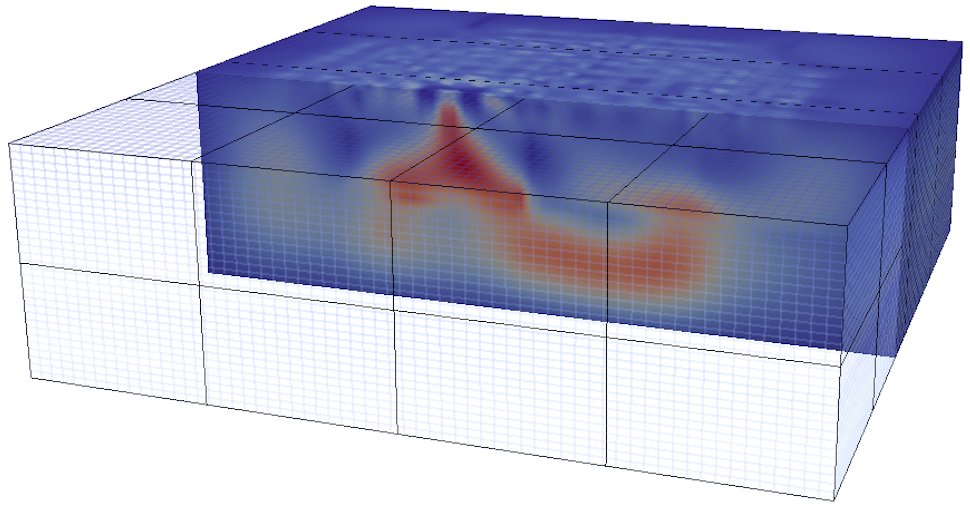}
       &
      \includegraphics[width=0.45\textwidth, height=1.6in]{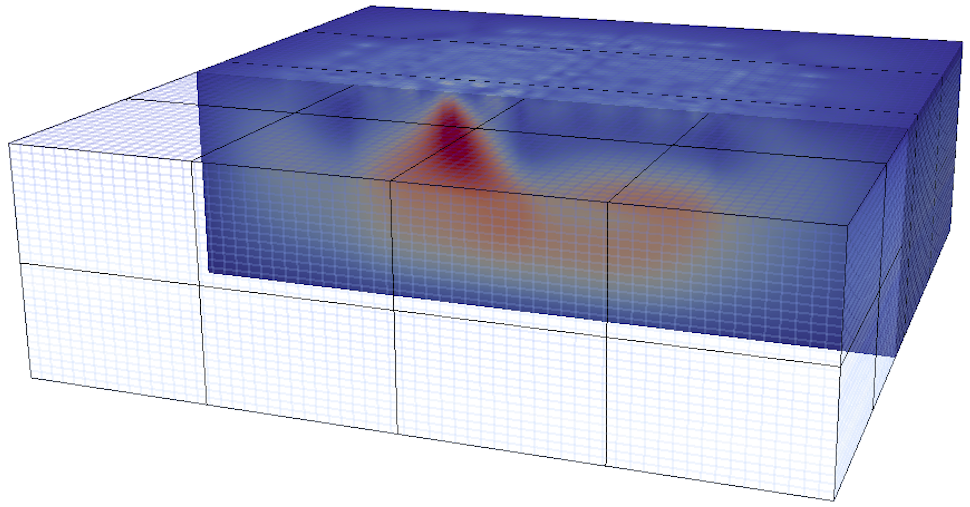} 
      \\ 
      \multicolumn{3}{c}{\includegraphics[width=0.7\textwidth, height=0.4in]{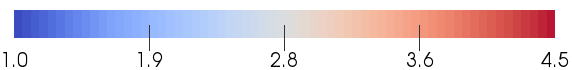}}
      \\ 
    \end{tabular}
      \caption{\textit{Model reconstructions for the SEG/EAGE test problem using the full mesh, and the adaptive and fixed multiscale inversions for two different coarsenings: $8^3$ and $16^3$ fine mesh cells per coarse cell. The figures were reproduced using the data visualization software Paraview \cite{ayachit2015paraview}.}}
  \label{fig:SEGModels}
\end{figure}

As a more realistic test problem, we consider the 3D SEG/EAGE model problem on a mesh of size $64 \times 64 \times 32$ with $72$ sources and $3698$ receivers. The test data is generated in \texttt{jInv} using the fine mesh. We use Lagrange polynomials to construct the multiscale basis and augment the basis by adding 72 basis functions of the global skeleton. We also use 862 and 3938 local basis boundary conditions described in Sec.~\ref{sub:MultiscaleFEM} for the coarse meshes consisting of $16^3$ and $8^3$ blocks, respectively. The construction of these boundary conditions required 72 fine mesh PDE solves in an offline phase which took about $2.21$ seconds. We compare the performance for a very coarse mesh containing blocks with $16^3$ fine mesh cells each and a more moderate coarsening using blocks of size $8^3$. For this test problem, we have $k=1009$ and $4415$ basis functions for the $16^3$ and $8^3$ coarsenings, respectively. Similar to the block model test problem, we use MUMPS and the block CG method with at most 100 iterations and stopping tolerance of $10^{-6}$ to solve the fine mesh forward problem, and we use MUMPS for the reduced multiscale forward problems. Again, the block CG stopping criteria are chosen such that the relative error of the block CG reconstruction has the same order of magnitude as the relative error of the multiscale inversions; see Table \ref{tab:SEGmodels}. For the inversions, we use 10 projected Gauss-Newton iterations with at most $15$ CG iterations for each step. We add $1\%$ noise to the data, and enforce smoothness by using a diffusion regularizer with regularization parameter $\alpha=10^{-15}$ in both full mesh and multiscale inversions.

In Table~\ref{tab:SEGmodels}, we show results for the SEG model test problem. The full mesh reconstruction requires solving a $139,425 \times 139,425$ linear system for the forward problem, whereas in the coarse meshes consisting of $8^3$ and $16^3$ cells, we project the PDEs down to $4415 \times 4415$ and $1009 \times 1009$ linear systems, respectively. As in the previous experiment, the runtimes of the adaptive inversions are larger than the one for MUMPS but considerably smaller than the one obtained using the iterative block CG algorithm as a linear solver. 

From the relative errors and Fig.~\ref{fig:SEGModels}, we again obtain more accurate reconstructions from the adaptive inversion - this is particularly clear in the $16^3$ reconstructions. In the $8^3$ case, the coarsening is fine enough so that the models are very similar; however, the adaptive multiscale reconstruction captures the peak of the full mesh model, whereas in the fixed multiscale reconstruction, there are discontinuities near the peak of the model. This is reflected in the relative errors in Table~\ref{tab:SEGmodels}.
\begin{table*}[h]
\centering
\caption{\textit{Strong scaling tests for constructing $\bfS_k, \bfY_k(\bfv, \bfm), \bfX_k(\bfw, \bfm)$ using a coarsening of $16^3$ fine cells per coarse cell, and with 72 sources and 3698 receivers. Computations are run on a Microway system that has four Intel Xeon E5-4627 CPUs with 40 cores and 1 TB of memory.}}
\label{tab:scalingTests}
\begin{tabular}{|c|c|c|c|c|c|c|}
\hline
& \multicolumn{3}{c|}{time (seconds)}                                      & \multicolumn{3}{c|}{speedup}                        \\ \hline
number of workers  & $\bfS_k(\bfm)$   & $\bfY_k(\bfv, \bfm)$  & $\bfX_k(\bfw, \bfm) $ & $\bfS_k $& $\bfY_k(\bfv, \bfm) $& $\bfX_k(\bfw, \bfm) $\\ \hline
1                    & 3.8595       & 2.7380                & 2.9933              & 1       & 1                   & 1                   \\ 
2                    & 2.7331       & 1.6775              & 1.8633                & 1.41    & 1.63                & 1.61                \\ 
4                    & 1.8164       & 1.0514              & 1.3548              & 2.12    & 2.60                & 2.21                \\ 
8                    & 1.4150       & 1.0658                & 1.3818                & 2.73    & 2.57                & 2.16                \\ \hline
\end{tabular}
\end{table*}
\subsection{Parallel Efficiency}
The computations of $\bfS_k(\bfm), \bfY_k(\bfv, \bfm), \bfX_k(\bfw, \bfm)$ require local PDE solves which can be performed independently on each coarse cell and provides an opportunity for parallel processing; see also Sec.~\ref{sec4}. We test the behavior of our implementation of these methods for a fixed problem size as we increase the number of workers from 1 to 8 as shown in Table~\ref{tab:scalingTests}. We use the same setup as in the example  from Sec.~\ref{sub:SEG}, i.e., a mesh of size $64 \times 64 \times 32$ with a coarsening of $16^3$ fine mesh cells per coarse cell. Increasing the number of workers from 1 to 8 decreases the runtimes of $\bfS_k(\bfm)$, $\bfY_k(\bfv, \bfm)$, and $\bfX_k(\bfw, \bfm)$ from around $3.86, 2.74$, and $2.99$ seconds to $1.41, 1.06$, and $1.38$ seconds, respectively. We therefore get speedup factors of $2.73$ for $\bfS_k(\bfm)$, $2.57$ for $\bfY_k(\bfv, \bfm)$, and $2.16$ for $\bfX_k(\bfw, \bfm)$ on the given machine, on which some resources such as caches are shared among workers. A summary of the runtimes can be found in Table~\ref{tab:scalingTests}. 

\section{Conclusion}
\label{sec6}
We embed a Multiscale Finite Volume (MSFV) methods into a PDE-constrained optimization framework and demonstrate its potential for solving high-dimensional parameter estimation problems. 
As usual in Model Order Reduction (MOR) techniques, we reduce the computational costs associated with the PDE constraint by projecting the discrete PDEs onto a lower-dimen\-sional subspace.
Following the MSFV approach, we obtain a reduced version of the original fine mesh PDE problem by projecting it onto a nested coarse mesh using an operator-dependent Galerkin projection. 
The key novelty of our method is the combination of MSFV and the numerical optimization scheme used for parameter estimation. Here,  we exploit the fact that the multiscale basis is sensitive to the current PDE parameter and propose a reduced inverse problem featuring an adaptive projection that can be solved using derivative-based optimization. We outline the potential of our method using two inverse conductivity problems in 3D that are inspired by Direct Current Resistivity. 

What sets our approach apart from existing works on ROM for PDE constrained optimization is the choice of multiscale methods for obtaining the reduced problem. This choice is mainly motivated by two reasons. First, using multiscale methods as ROM techniques avoids the necessity of sampling the parameter space, which is problematic in high-dimensional spaces. Apart from optional boundary conditions, the method is fully online and does not require solving the fine mesh problem. Second, multiscale methods simplify the formulation and solution of adaptive inversion  problems as the multiscale basis vary smoothly with respect to the PDE parameters. Similar to the recent work of \cite{de2017multiscale}, we show that computing derivatives is tractable by using local sensitivity computations; see Sec.~\ref{sec4}. 

Following a discretize-then-optimize strategy, we explicitly differentiate the solution of the discretized problem obtained with the MSFV solver with respect to the parameter to be estimated. 
Earlier works on differentiating multiscale solutions obtained approximate derivatives through interpolation~\cite{krogstad2011adjoint}, used adjoint-based approaches~\cite{fu2009multiscale, fu2010multiscale,fu2011multiscale}, or required automatic differentiation~\cite{de2017multiscale}. 
Similar to other discretize-then-optimize approaches, our method produces an accurate gradient of the discrete objective problem irrespective of the mesh size and quality of the reduced order model. We also demonstrate that the involved local sensitivity computations can be performed in parallel providing additional opportunity for speedup.

Our numerical experiments show that the adaptive multiscale provides parameter estimates that feature details below the coarse mesh resolution. As expected (see discussion in Sec.~\ref{sub:optimization_with_fixed_reduced_space}), our numerical experiments confirm that ignoring the dependence of $\bfS_k$ on the model parameter $\bfm$ also degrades reconstruction quality, especially when using aggressive coarsening.
The computational complexity of the multiscale inversion grows linearly with respect to the number of coarse mesh blocks and (if direct methods are used) cubic with respect to the number of fine mesh cells in each coarse mesh blocks. Thus, larger computational savings are to be expected as the number of fine mesh cells grows and in particular for problems where direct linear solvers cannot be applied and iterative solvers are required. This can also be seen in our numerical experiments. 

We applied our method to the DCR survey, where the forward problem involves solving the diffusion equation. We intend to further explore this method for other parameter estimation problems arising from electromagnetic, and gravity based surveys which involve solving different PDEs \cite{caudillo2017framework,wilhelms20fast}. 
Multiscale Finite Volume (MSFV)  methods provide more flexibility in building the reduced basis than used in this work. For example, additional basis functions representing sources and receivers can be added; see Sec.~\ref{sub:MultiscaleFEM}. Another approach to reduce the impact of boundary conditions on the multiscale basis is using oversampling~\cite{ElfendievHou2009,caudillo2016oversampling}. 

\section{Acknowledgements}
We would like to thank Eldad Haber (University of British Columbia, Department of Earth, Ocean and Atmospheric Sciences) for his insight and helpful suggestions. This work is supported by National Science Foundation (NSF) award DMS 1522599.

\bibliographystyle{abbrv}      
\bibliography{MSPaper.bib}

\begin{thebibliography}{10}

\bibitem{amestoy2000mumps}
P.~R. Amestoy, I.~S. Duff, J.-Y. L’Excellent, and J.~Koster.
\newblock Mumps: a general purpose distributed memory sparse solver.
\newblock In {\em International Workshop on Applied Parallel Computing}, pages
  121--130. Springer, 2000.

\bibitem{aminzadeh1997}
F.~Aminzadeh, B.~Jean, and T.~Kunz.
\newblock {\em 3-D salt and overthrust models}.
\newblock Society of Exploration Geophysicists, 1997.

\bibitem{AmsallemFarhat2011}
D.~Amsallem and C.~Farhat.
\newblock {An online method for interpolating linear parametric reduced-order
  models}.
\newblock {\em SIAM Journal on Scientific Computing}, 2011.

\bibitem{Arridge1999}
S.~R. Arridge.
\newblock {Optical tomography in medical imaging}.
\newblock {\em Inverse Problems}, 15(2):R41--R93, 1999.

\bibitem{ArridgeSchotland2009}
S.~R. Arridge and J.~C. Schotland.
\newblock {Optical tomography: forward and inverse problems}.
\newblock {\em Inverse Problems}, 25(12):123010--60, 2009.

\bibitem{ayachit2015paraview}
U.~Ayachit.
\newblock The paraview guide: a parallel visualization application.
\newblock 2015.

\bibitem{bai2005reduced}
Z.~Bai, P.~M. Dewilde, and R.~W. Freund.
\newblock Reduced-order modeling.
\newblock {\em Handbook of numerical analysis}, 13:825--895, 2005.

\bibitem{BarraultEtAl2004}
M.~Barrault, Y.~Maday, N.~C. Nguyen, and A.~T. Patera.
\newblock {An {\textquoteleft}empirical interpolation{\textquoteright} method:
  application to efficient reduced-basis discretization of partial differential
  equations}.
\newblock {\em Comptes Rendus Mathematique}, 339(9):667--672, 2004.

\bibitem{BaurEtAl2011}
U.~Baur, C.~Beattie, P.~Benner, and S.~Gugercin.
\newblock {Interpolatory Projection Methods for Parameterized Model Reduction}.
\newblock {\em SIAM Journal on Scientific Computing}, 33(5):2489--2518, Jan.
  2011.

\bibitem{bezanson2017julia}
J.~Bezanson, A.~Edelman, S.~Karpinski, and V.~B. Shah.
\newblock Julia: A fresh approach to numerical computing.
\newblock {\em SIAM Review}, 59(1):65--98, 2017.

\bibitem{schulz2011computational}
A.~Borzi and V.~Schulz.
\newblock {\em Computational optimization of systems governed by partial
  differential equations}.
\newblock SIAM, Philadelphia, 2011.

\bibitem{BuiThanhEtAl2008}
T.~Bui-Thanh, K.~Willcox, and O.~Ghattas.
\newblock {Model Reduction for Large-Scale Systems with High-Dimensional
  Parametric Input Space}.
\newblock {\em SIAM Journal on Scientific Computing}, 30(6):3270--3288, Jan.
  2008.

\bibitem{calo2016randomized}
V.~M. Calo, Y.~Efendiev, J.~Galvis, and G.~Li.
\newblock Randomized oversampling for generalized multiscale finite element
  methods.
\newblock {\em Multiscale Modeling \& Simulation}, 14(1):482--501, 2016.

\bibitem{caudillo2016oversampling}
L.~Caudillo~Mata, E.~Haber, and C.~Schwarzbach.
\newblock An oversampling technique for multiscale finite volume method to
  simulate frequency-domain electromagnetic responses.
\newblock In {\em SEG Technical Program Expanded Abstracts 2016}, pages
  981--985. Society of Exploration Geophysicists, 2016.

\bibitem{caudillo2017framework}
L.~A. Caudillo-Mata, E.~Haber, L.~J. Heagy, and C.~Schwarzbach.
\newblock A framework for the upscaling of the electrical conductivity in the
  quasi-static maxwell’s equations.
\newblock {\em Journal of Computational and Applied Mathematics}, 317:388--402,
  2017.

\bibitem{CheneyEtAl1999}
M.~Cheney, D.~Isaacson, and J.~C. Newell.
\newblock {Electrical impedance tomography}.
\newblock {\em SIAM review}, 41(1):85--101, 1999.

\bibitem{chung2016adaptive}
E.~Chung, Y.~Efendiev, and T.~Y. Hou.
\newblock Adaptive multiscale model reduction with generalized multiscale
  finite element methods.
\newblock {\em Journal of Computational Physics}, 320:69--95, 2016.

\bibitem{chung2003convergence}
E.~T. Chung, Q.~Du, and J.~Zou.
\newblock Convergence analysis of a finite volume method for maxwell's
  equations in nonhomogeneous media.
\newblock {\em SIAM Journal on Numerical Analysis}, 41(1):37--63, 2003.

\bibitem{chung2011energy}
E.~T. Chung, Y.~Efendiev, and R.~L. Gibson~Jr.
\newblock An energy-conserving discontinuous multiscale finite element method
  for the wave equation in heterogeneous media.
\newblock {\em Advances in Adaptive Data Analysis}, 3(01n02):251--268, 2011.

\bibitem{chung2015mixed}
E.~T. Chung, Y.~Efendiev, and C.~S. Lee.
\newblock Mixed generalized multiscale finite element methods and applications.
\newblock {\em Multiscale Modeling \& Simulation}, 13(1):338--366, 2015.

\bibitem{chung2016generalized}
E.~T. Chung, Y.~Efendiev, G.~Li, and M.~Vasilyeva.
\newblock Generalized multiscale finite element methods for problems in
  perforated heterogeneous domains.
\newblock {\em Applicable Analysis}, 95(10):2254--2279, 2016.

\bibitem{de2017multiscale}
R.~J. de~Moraes, J.~R. Rodrigues, H.~Hajibeygi, and J.~D. Jansen.
\newblock Multiscale gradient computation for flow in heterogeneous porous
  media.
\newblock {\em Journal of Computational Physics}, 336:644--663, 2017.

\bibitem{deSturlerEtAl2013}
E.~de~Sturler, S.~Gugercin, M.~E. Kilmer, S.~Chaturantabut, C.~Beattie, and
  M.~O'Connell.
\newblock {Nonlinear Parametric Inversion using Interpolatory Model Reduction}.
\newblock {\em arXiv.org}, Nov. 2013.

\bibitem{de2015nonlinear}
E.~De~Sturler, S.~Gugercin, M.~E. Kilmer, S.~Chaturantabut, C.~Beattie, and
  M.~O'Connell.
\newblock Nonlinear parametric inversion using interpolatory model reduction.
\newblock {\em SIAM Journal on Scientific Computing}, 37(3):B495--B517, 2015.

\bibitem{OliverBook2008}
O.~S. Dean, A.~C. Reynolds, and N.~Liu.
\newblock {\em Inverse Theory for Petroleum Reservoir Characterization and
  History Matching}.
\newblock Cambridge University Press, Cambridge, 2008.

\bibitem{deymor}
A.~Dey and H.~Morrison.
\newblock Resistivity modeling for arbitrarily shaped three dimensional
  structures.
\newblock {\em Geophysics}, 44(4):753--780, 1979.

\bibitem{durlofsky2007adaptive}
L.~Durlofsky, Y.~Efendiev, and V.~Ginting.
\newblock An adaptive local--global multiscale finite volume element method for
  two-phase flow simulations.
\newblock {\em Advances in Water Resources}, 30(3):576--588, 2007.

\bibitem{ElfendievHou2009}
Y.~Efendiev and T.~Y. Hou.
\newblock {\em Multiscale finite element methods: theory and applications},
  volume~4.
\newblock Springer Science \& Business Media, 2009.

\bibitem{efendiev2000convergence}
Y.~R. Efendiev, T.~Y. Hou, and X.-H. Wu.
\newblock Convergence of a nonconforming multiscale finite element method.
\newblock {\em SIAM Journal on Numerical Analysis}, 37(3):888--910, 2000.

\bibitem{ElmanLiao2013}
H.~C. Elman and Q.~Liao.
\newblock {Reduced Basis Collocation Methods for Partial Differential Equations
  with Random Coefficients}.
\newblock {\em SIAM/ASA Journal on Uncertainty Quantification}, 1(1):192--217,
  Jan. 2013.

\bibitem{EpanomeritakisAkcelikGhattasBielak2008}
I.~Epanomeritakis, V.~Akcelik, O.~Ghattas, and J.~Bielak.
\newblock A newton-cg method for large-scale three-dimensional elastic
  full-waveform seismic inversion.
\newblock {\em Inverse Problems}, 24(3):034015, 2008.

\bibitem{FahlEtAl2000}
M.~Fahl and E.~W. Sachs.
\newblock {Trust-region Proper Orthogonal Decomposition for Flow Control},
  2000.

\bibitem{FengBenner2007}
L.~Feng and P.~Benner.
\newblock {A robust algorithm for parametric model order reduction}.
\newblock {\em PAMM}, 7(1):1021501--1021502, Dec. 2007.

\bibitem{fu2011multiscale}
J.~Fu, J.~Caers, and H.~A. Tchelepi.
\newblock A multiscale method for subsurface inverse modeling: Single-phase
  transient flow.
\newblock {\em Advances in water resources}, 34(8):967--979, 2011.

\bibitem{fu2009multiscale}
J.~Fu, H.~Tchelepi, and J.~Caers.
\newblock A multiscale adjoint method for large-scale sensitivity computation
  in subsurface flow simulation, 2009.

\bibitem{fu2010multiscale}
J.~Fu, H.~A. Tchelepi, and J.~Caers.
\newblock A multiscale adjoint method to compute sensitivity coefficients for
  flow in heterogeneous porous media.
\newblock {\em Advances in water resources}, 33(6):698--709, 2010.

\bibitem{GalballyEtAl2009}
D.~Galbally, K.~Fidkowski, K.~Willcox, and O.~Ghattas.
\newblock {Non-linear model reduction for uncertainty quantification in
  large-scale inverse problems}.
\newblock {\em International Journal for Numerical Methods in Engineering},
  81:1581--1608, 2009.

\bibitem{gallivan2004model}
K.~Gallivan, A.~Vandendorpe, and P.~Van~Dooren.
\newblock Model reduction of mimo systems via tangential interpolation.
\newblock {\em SIAM Journal on Matrix Analysis and Applications},
  26(2):328--349, 2004.

\bibitem{gamkrelidze2013principles}
R.~Gamkrelidze.
\newblock {\em Principles of optimal control theory}, volume~7.
\newblock Springer Science \& Business Media, 2013.

\bibitem{ghasemi2015fast}
M.~Ghasemi, Y.~Yang, E.~Gildin, Y.~R. Efendiev, and V.~M. Calo.
\newblock Fast multiscale reservoir simulations using pod-deim model reduction.
\newblock In {\em SPE reservoir simulation symposium}, 2015.

\bibitem{grimme1996model}
E.~J. Grimme, D.~C. Sorensen, and P.~Van~Dooren.
\newblock Model reduction of state space systems via an implicitly restarted
  lanczos method.
\newblock {\em Numerical algorithms}, 12(1):1--31, 1996.

\bibitem{HaasdonkOhlberger2008}
B.~Haasdonk and M.~Ohlberger.
\newblock {Reduced basis method for finite volume approximations of
  parametrized linear evolution equations}.
\newblock {\em ESAIM: Mathematical Modelling and Numerical Analysis},
  42(02):277--302, Mar. 2008.

\bibitem{haber2014computational}
E.~Haber.
\newblock {\em Computational methods in geophysical electromagnetics}.
\newblock SIAM, 2014.

\bibitem{haber2016solving}
E.~Haber.
\newblock Solving csem problems with massive number of sources and receivers.
\newblock In {\em 78th EAGE Conference and Exhibition 2016}, 2016.

\bibitem{haber2012adaptive}
E.~Haber, E.~Holtham, J.~Granek, D.~Marchant, D.~Oldenburg, C.~Schwarzbach, and
  R.~Shekhtman.
\newblock An adaptive mesh method for electromagnetic inverse problems.
\newblock In {\em SEG Technical Program Expanded Abstracts 2012}, pages 1--6.
  Society of Exploration Geophysicists, 2012.

\bibitem{haber2014multiscale}
E.~Haber and L.~Ruthotto.
\newblock A multiscale finite volume method for maxwell's equations at low
  frequencies.
\newblock {\em Geophysical Journal International}, 199(2):1268--1277, 2014.

\bibitem{hajibeygi2011adaptive}
H.~Hajibeygi and P.~Jenny.
\newblock Adaptive iterative multiscale finite volume method.
\newblock {\em Journal of Computational Physics}, 230(3):628--643, 2011.

\bibitem{hajibeygi2014compositional}
H.~Hajibeygi, H.~A. Tchelepi, et~al.
\newblock Compositional multiscale finite-volume formulation.
\newblock {\em SPE Journal}, 19(02):316--326, 2014.

\bibitem{HimpeOhlberger2015}
C.~Himpe and M.~Ohlberger.
\newblock {Data-driven combined state and parameter reduction for inverse
  problems}.
\newblock {\em Advances in Computational Mathematics}, 41(5):1343--1364, 2015.

\bibitem{HoreshHaber2011}
L.~Horesh and E.~Haber.
\newblock {A Second Order Discretization of Maxwell's Equations in the
  Quasi-Static Regime on OcTree Grids}.
\newblock {\em SIAM Journal on Scientific Computing}, 33(5):2805–--2822,
  2011.

\bibitem{hou1999convergence}
T.~Hou, X.-H. Wu, and Z.~Cai.
\newblock Convergence of a multiscale finite element method for elliptic
  problems with rapidly oscillating coefficients.
\newblock {\em Mathematics of Computation of the American Mathematical
  Society}, 68(227):913--943, 1999.

\bibitem{HouWu1997}
T.~Y. Hou and X.~H. Wu.
\newblock {A multiscale finite element method for elliptic problems in
  composite materials and porous media}.
\newblock {\em Journal of Computational Physics}, 134(1):169--189, 1997.

\bibitem{jenny2003multi}
P.~Jenny, S.~Lee, and H.~A. Tchelepi.
\newblock Multi-scale finite-volume method for elliptic problems in subsurface
  flow simulation.
\newblock {\em Journal of Computational Physics}, 187(1):47--67, 2003.

\bibitem{jenny2005adaptive}
P.~Jenny, S.~H. Lee, and H.~A. Tchelepi.
\newblock Adaptive multiscale finite-volume method for multiphase flow and
  transport in porous media.
\newblock {\em Multiscale Modeling \& Simulation}, 3(1):50--64, 2005.

\bibitem{jenny2009modeling}
P.~Jenny and I.~Lunati.
\newblock Modeling complex wells with the multi-scale finite-volume method.
\newblock {\em Journal of Computational Physics}, 228(3):687--702, 2009.

\bibitem{jiang2007multiscale}
L.~Jiang, Y.~Efendiev, and V.~Ginting.
\newblock Multiscale methods for parabolic equations with continuum spatial
  scales.
\newblock {\em Discrete and Continuous Dynamical Systems Series B}, 8(4):833,
  2007.

\bibitem{KaipioSomersalo2006}
J.~Kaipio and E.~Somersalo.
\newblock {\em {Statistical and Computational Inverse Problems}}, volume 160 of
  {\em Applied Mathematical Sciences}.
\newblock Springer Science {\&} Business Media, New York, Mar. 2006.

\bibitem{KalchevEtAl2016}
D.~Z. Kalchev, C.~S. Lee, U.~Villa, and Y.~Efendiev.
\newblock {Upscaling of mixed finite element discretization problems by the
  spectral AMGe method}.
\newblock {\em SIAM Journal on {\ldots}}, 38(5):A2912--A2933, 2016.

\bibitem{KaletaEtAl2010}
M.~P. Kaleta, R.~G. Hanea, A.~W. Heemink, and J.-D. Jansen.
\newblock {Model-reduced gradient-based history matching}.
\newblock {\em Computational Geosciences}, 15(1):135--153, Aug. 2010.

\bibitem{krogstad2011adjoint}
S.~Krogstad, V.~L. Hauge, A.~Gulbransen, et~al.
\newblock Adjoint multiscale mixed finite elements.
\newblock {\em SPE Journal}, 16(01):162--171, 2011.

\bibitem{KunischVolkwein2008}
K.~Kunisch and S.~Volkwein.
\newblock {Proper orthogonal decomposition for optimality systems}.
\newblock {\em ESAIM: Mathematical Modelling and Numerical Analysis},
  42(1):1--23, 2008.

\bibitem{lee2009adaptive}
S.~H. Lee, H.~Zhou, and H.~A. Tchelepi.
\newblock Adaptive multiscale finite-volume method for nonlinear multiphase
  transport in heterogeneous formations.
\newblock {\em Journal of Computational Physics}, 228(24):9036--9058, 2009.

\bibitem{LiebermanEtAl2010}
C.~Lieberman, K.~Willcox, and O.~Ghattas.
\newblock {Parameter and State Model Reduction for Large-Scale Statistical
  Inverse Problems}.
\newblock {\em SIAM Journal on Scientific Computing}, 32(5):2523--2542, Jan.
  2010.

\bibitem{LipnikovEtAl2004}
K.~Lipnikov, J.~Morel, and M.~Shashkov.
\newblock {Mimetic finite difference methods for diffusion equations on
  non-orthogonal non-conformal meshes}.
\newblock {\em Journal of Computational Physics}, 199(2):589--597, 2004.

\bibitem{LipponenEtAl2013}
A.~Lipponen, A.~Sepp{\"a}nen, and J.~Kaipio.
\newblock {Electrical impedance tomography imaging with reduced-order model
  based on proper orthogonal decomposition}.
\newblock {\em Journal of Electronic Imaging}, 22(2):023008--16, 2013.

\bibitem{LohmannEid2007}
B.~Lohmann and R.~Eid.
\newblock {Efficient order reduction of parametric and nonlinear models by
  superposition of locally reduced models}.
\newblock {\em Methoden und Anwendungen der Regelungstechnik}, 2007.

\bibitem{lunati2008multiscale}
I.~Lunati and P.~Jenny.
\newblock Multiscale finite-volume method for density-driven flow in porous
  media.
\newblock {\em Computational Geosciences}, 12(3):337--350, 2008.

\bibitem{MacLachlanMoulton2006}
S.~P. MacLachlan and J.~D. Moulton.
\newblock {Multilevel upscaling through variational coarsening}.
\newblock {\em Water Resources Research}, 42(2):131--9, 2006.

\bibitem{MartinEtAl2012}
J.~Martin, L.~C. Wilcox, C.~Burstedde, and O.~Ghattas.
\newblock {A Stochastic Newton MCMC Method for Large-Scale Statistical Inverse
  Problems with Application to Seismic Inversion}.
\newblock {\em SIAM Journal on Scientific Computing}, 34(3):A1460--A1487, 2012.

\bibitem{mcgillivray1992forward}
P.~R. McGillivray.
\newblock {\em Forward modeling and inversion of DC resistivity and MMR data}.
\newblock PhD thesis, University of British Columbia, 1992.

\bibitem{moyner2016multiscale}
O.~M{\o}yner and K.-A. Lie.
\newblock A multiscale restriction-smoothed basis method for high contrast
  porous media represented on unstructured grids.
\newblock {\em Journal of Computational Physics}, 304:46--71, 2016.

\bibitem{NegriEtAl2013}
F.~Negri, G.~Rozza, A.~Manzoni, and A.~Quarteroni.
\newblock {Reduced basis method for parametrized elliptic optimal control
  problems}.
\newblock {\em SIAM Journal on Scientific Computing}, 35(5):A2316--A2340, 2013.

\bibitem{OConnellEtAl2017}
M.~O'Connell, M.~E. Kilmer, and E.~de~Sturler.
\newblock {Computing Reduced Order Models via Inner-Outer Krylov Recycling in
  Diffuse Optical Tomography}.
\newblock {\em SIAM Journal on Scientific Computing}, 39(2):B272--B297, 2017.

\bibitem{o1980block}
D.~P. O'Leary.
\newblock The block conjugate gradient algorithm and related methods.
\newblock {\em Linear algebra and its applications}, 29:293--322, 1980.

\bibitem{PanzerEtAl2010}
H.~Panzer, J.~Mohring, R.~Eid, and B.~Lohmann.
\newblock {Parametric Model Order Reduction by Matrix Interpolation}.
\newblock {\em at - Automatisierungstechnik}, 58(8), 2010.

\bibitem{Parker1994}
R.~L. Parker.
\newblock {\em Geophysical Inverse Theory}.
\newblock Princeton University Press, Princeton NJ, 1994.

\bibitem{parramore2016multiscale}
E.~Parramore, M.~G. Edwards, M.~Pal, and S.~Lamine.
\newblock Multiscale finite-volume cvd-mpfa formulations on structured and
  unstructured grids.
\newblock {\em Multiscale Modeling \& Simulation}, 14(2):559--594, 2016.

\bibitem{BennerEtAl2014}
S.~V. Peter~Brenner, Ekkehard~Sachs.
\newblock {Model Order Reduction for PDE Constrained Optimization}.
\newblock pages 1--25, 2014.

\bibitem{pratt1999}
R.~Pratt.
\newblock Seismic waveform inversion in the frequency domain, part 1: Theory,
  and verification in a physical scale model.
\newblock {\em Geophysics}, 64(3):888--901, 1999.

\bibitem{jInvRuthottoHaberTreister}
L.~Ruthotto, E.~Treister, and E.~Haber.
\newblock jinv - a flexible julia package for {PDE} parameter estimation.
\newblock {\em CoRR}, abs/1606.07399, 2016.

\bibitem{saad2003iterative}
Y.~Saad.
\newblock {\em Iterative methods for sparse linear systems}.
\newblock SIAM, 2003.

\bibitem{sargent2000optimal}
R.~Sargent.
\newblock Optimal control.
\newblock {\em Journal of Computational and Applied Mathematics},
  124(1):361--371, 2000.

\bibitem{SpantiniEtAl2015}
A.~Spantini, A.~Solonen, T.~Cui, J.~Martin, L.~Tenorio, and Y.~Marzouk.
\newblock {Optimal Low-Rank Approximations of Bayesian Linear Inverse
  Problems}.
\newblock {\em SIAM Journal on Scientific Computing}, 37(6):A2451--A2487, 2015.

\bibitem{van2015penalty}
T.~van Leeuwen and F.~J. Herrmann.
\newblock A penalty method for pde-constrained optimization in inverse
  problems.
\newblock {\em Inverse Problems}, 32(1):015007, 2015.

\bibitem{Wardhow1988}
S.~Ward and G.~Hohmann.
\newblock Electromagnetic theory for geophysical applications.
\newblock {\em Electromagnetic Methods in Applied Geophysics}, 1:131--311,
  1988.
\newblock Soc. Expl. Geophys.

\bibitem{wilhelms20fast}
W.~Wilhelms, C.~Schwarzbach, R.-U. B{\"o}rner, and K.~Spitzer.
\newblock A fast 3d mt inversion--the forward operator behind.
\newblock {\em Journal of Future Generation Computer Systems}, 20(3):475--487.

\bibitem{WillcoxPeraire2002}
K.~Willcox and J.~Peraire.
\newblock {Balanced Model Reduction via the Proper Orthogonal Decomposition}.
\newblock {\em AIAA Journal}, 40(11):2323--2330, Nov. 2002.

\end{thebibliography}


\end{document}